\documentclass[a4paper,10pt]{article}
\usepackage{amsmath,amssymb}
\usepackage[latin1]{inputenc}  
\begin{document}
\title{\bf Smoothing properties for the higher order nonlinear Schr\"{o}dinger equation with
constant coefficients}
\author{Mauricio Sep\'{u}lveda \thanks{
Departamento de Ingenier\'{\i}a Matem\a'{a}tica, Universidad de
Concepci\a'{o}n, Casilla 160-C, Concepci\a'{o}n, Chile.
mauricio@ing-mat.udec.cl}\quad \quad Octavio Vera
Villagr\a'{a}n.\thanks{Departamento de Matem\a'{a}tica,
Universidad del B\a'{\i}o-B\a'{\i}o, Collao 1202, Casilla 5-C,
Concepci\a'{o}n, Chile. overa@ubiobio.cl}}
\date{}
\maketitle
\begin{abstract}
\noindent We study local and  global existence and smoothing
properties for the initial value problem associated to a higher
order nonlinear Schr\"{o}dinger equation with constant coefficients
which appears as a model for propagation of pulse in optical fiber.
\end{abstract}
\noindent \underline{Keywords and phrases}:  Evolution equations,
weighted Sobolev space, gain in regularity.\\
\\
\noindent Mathematics Subject Classification: {35Q53,\,47J35}

\renewcommand{\theequation}{\thesection.\arabic{equation}}
\setcounter{equation}{0}\section{Introduction}
We consider the initial value problem
\[(P) \left \{ \begin{array}{ll}
i\,u_{t} + \omega \,u_{xx} + i\,\beta \,u_{xxx} + |u|^{2}\,u=0 &
\quad x,\,t\in \mathbb{R}\\
u(x,\,0) = u_{0}(x) &
\end{array}
\right. \] where $\omega ,\,\beta \in \mathbb{R},$ $\beta \neq 0 $
and $u=u(x,t)$ is a complex valued function. The above equation is a
particular case of the equation
\[(Q) \left \{ \begin{array}{ll}
i\,u_{t} + \omega \,u_{xx} + i\,\beta \,u_{xxx} + \gamma \,|u|^{2}\,u +
i\,\delta \,|u|^{2}\,u_{x} + i\,\epsilon \,u^{2}\,\overline{u}_{x}=0 &
\quad x,\,t\in \mathbb{R}\\
u(x,\,0)  = u_{0}(x) &
\end{array}
\right. \] where $\omega ,\,\beta ,\,\gamma ,\,\delta $ are real
numbers with $\beta \neq 0.$ This equation was first proposed by A.
Hasegawa and Y. Kodama \cite{ha1} as a model for the propagation of
a signal in an optic fiber (see also \cite{ko1}). The equation $(Q)$
can be reduced to other well known equations. For instance, setting
$\omega =1,$ $\beta = \delta =\epsilon =0$ in $(Q)$ we have the
semilinear Schr\"{o}dinger equation, i. e.,
\begin{eqnarray*}
i\,u_{t} + u_{xx} + \gamma \,|u|^{2}\,u =0.\qquad (Q_{1})
\end{eqnarray*}
If we let $\beta = \gamma =0$ and $\omega =1$ in $(Q),$ we obtain
the derivative nonlinear Schr\"{o}dinger equation
\begin{eqnarray*}
i\,u_{t} + u_{xx} + i\,\delta \,|u|^{2}\,u_{x} + i\,\epsilon
\,u^{2}\,\overline{u}_{x}=0.\qquad (Q_{2})
\end{eqnarray*}
Letting $\alpha = \gamma = \epsilon =0$ in $(Q),$ the equation that arises is the
complex modified Korteweg-de Vries equation,
\begin{eqnarray*}
i\,u_{t} + i\,\beta \,u_{xxx} + i\,\delta \,|u|^{2}\,u_{x} =0.\qquad
(Q_{3})
\end{eqnarray*}
The initial value problem for the equations $(Q_{1}),$ $(Q_{2})$ and
$(Q_{3})$ has been extensively studied in the last few years. See,
for instance, \cite{bi1,bi2,bo1,ca1,ca2,cr1,cr2,ka1,ke1,sa1,sj1} and
references therein. In 1992, C. Laurey \cite{la1} considered the
equation $(Q)$ and proved local well-posedness of the initial value
problem associated for data in $H^{s}(\mathbb{R}),$ $s>3/4,$ and
global well-posedness in $H^{s}(\mathbb{R}),$ $s\geq 1.$ In 1997, G.
Staffilani \cite{st1} for $(Q)$ established local well-posedness for
data in $H^{s}(\mathbb{R}),$ $s\geq 1/4$ improving Laurey's result.
A similar result was given in \cite{ca1,ca2}
with $w(t),$ $\beta(t)$ real functions.\\
Our aim in this paper, is to study gain in regularity for the
equation $(P).$ Specifically, we prove conditions on $(P)$ for which
initial data $u_{0}$ possessing sufficient decay at infinity and
minimal amount of regularity will lead to a unique solution $u(t)\in
C^{\infty}(\mathbb{R})$ for $0<t<T,$ where $T$ is the existence time
of the solution. We are not considering the equation $(Q)$ because
of the technique used here, we shall see that the last two terms in
$(Q)$ are not outstanding in the main inequality, indeed the two
last terms are
observed in the last two terms in the main inequality.\\
In 1986, N. Hayashi {\it et al.} \cite{ha1} showed that for the
nonlinear Schr\"{o}dinger equation (NLS): $i\,u_{t} + u_{xx} =
\lambda \,|u|^{p - 1}\,u,$ $(x,\,t)\in \mathbb{R} \times \mathbb{R}$
with initial condition $u(x,\,0)=u_{0}(x),$ $x\in \mathbb{R}$ and a
certain assumption on $\lambda $ and $p,$  all solutions of finite
energy are smooth for $t\neq 0$ provided the initial functions in
$H^{1}(\mathbb{R})$(or on $L^{2}(\mathbb{R})$) decay sufficiently
fast as $|x|\rightarrow \infty .$ The main tool is the operator $J$
defined by
$Ju=e^{i\,x^{2}/4\,t}\,(2\,i\,t)\,\partial_{x}(e^{-\,i\,x^{2}/4\,t}\,u)=
(x + 2\,i\,t\,\partial_{x})u$ which has the remarkable property that
it commutes with the operator $L$ defined by $L=(i\,\partial_{t} +
\partial_{x}^{2}),$ namely
$LJ - JL =[L,\,J]=0.$\\
For the Korteweg-de Vries type equation (KdV), J. C. Saut and M.
Temam \cite{sa1} remarked that a solution $u$ cannot gain or lose
regularity. They showed that if $u(x,\,0)=u_{0}(x)\in
H^{s}(\mathbb{R})$ for $s\geq 2,$ then $u(\,\cdot\,,\,t)\in
H^{s}(\mathbb{R})$ for all $t>0.$ For the KdV equation on the line,
Kato \cite{ka1} motivated by work of Cohen \cite{co1} showed that if
$u(x,\,0)=u_{0}(x)\in L_{b}^{2}\equiv H^{2}(\mathbb{R})\cap
L^{2}$($e^{bx}\,dx$)($b>0$) then the solution $u(x,\,t)$ of the KdV
equation becomes $C^{\infty}$ for all $t>0.$ A main ingredient in
the proof was the fact that formally the semi-group
$S(t)=e^{-\,\partial_{x}^{3}}$ in $L_{b}^{2}(\mathbb{R})$ is
equivalent to $S_{b}(t)= e^{-\,t\,(\partial_{x} - b)^{3}}$ in
$L^{2}(\mathbb{R})$ when $t>0.$ One would be inclined to believe
that this was a special property of the KdV equation. However, his
is not the case. The effect is due to the dispersive nature of the
linear part of the equation. Kruzkov and Faminskii \cite{kr1} proved
that $u(x,\,0)=u_{0}(x)\in L^{2}(\mathbb{R})$ such that
$x^{\alpha}\,u_{0}(x)\in L^{2}((0,\,+\infty)),$ the weak solution of
the KdV equation, has $l$-continuous space derivatives for all $t>0$
if $l<2\,\alpha .$ The proof of this result is based on the
asymptotic behavior of the Airy function and its derivatives, and on
the smoothing effect of the KdV equation which was found in
\cite{ka1,kr1}. While the proof of Kato appears to depend on special
a priori estimates, some of this mystery has been solved by the
result of local gain of finite regularity for various others linear
and nonlinear dispersive equations due to Ginibre and Velo
\cite{gi1} and others.
However, all of them require growth conditions on the nonlinear term. \\
In 1992, W. Craig, T. Kappeler and W. Strauss \cite{cr1,cr2} proved
for the fully nonlinear KdV equation $u_{t} + $
$f(u_{xxx},\,u_{xx},\,u_{x},\,u,\,x,\,t)=0,$ $x\in \mathbb{R},$
$t>0$ and certain additional assumption over $f$ that $C^{\infty}$
solutions $u(x,\,t)$ are obtained for all $t>0$ if the initial data
$u_{0}(x)$ decays faster than polynomially on $\mathbb{R}^{+}=\{x
\in \mathbb{R}:\;x>0\}$ and has certain initial Sobolev regularity.
Following this idea, H. Cai \cite{ca0} studied the nonlinear
equation of KdV-type of the form $u_{t} + u_{xxx} +
a(x,,t)\,f(u_{xx},\,u_{x},\,u,\,x,\,t)=0,$ where $a(x,\,t)$ is
positive and bounded, obtaining the same conclusion. Subsequent
works were given by O. Vera \cite{ve1,ve2,ve3,ve4} for a nonlinear
dispersive evolution equation, a KdV-Burgers type equation and for
KdV-Kawahara type equation, respectively. In more than one spatial
dimension, J. Levandosky \cite{le1}, proved infinite gain in
regularity results for nonlinear third-order equations. While
\cite{cr1} included local smoothing results for some mth-order
dispersive equation in $n$ spatial dimension, their results and the
techniques are different from those presented by Levandosky. First,
they consider equations with only a mild solution and Levandosky
considers equations with very general nonlinearities including a
fully nonlinear equation of the form
\begin{eqnarray*}
&  & u_{t} + f(D^{3}u,\,D^{2}u,\,Du,\,u,\,x,\,t)=0,\\
&  & u(x,\,y,\,0)=u_{0}(x,\,y).
\end{eqnarray*}
Secondly, they indicate local gain in finite regularity and
Levandosky proved complementary results showing the relationship
between the decay at infinity of the initial data and the amount
of gain in regularity. More specifically, it is proved a condition
under which an equation of the form
\begin{eqnarray*}
&  & u_{t} + a\,u_{xxx} + b\,u_{xxy} + c\,u_{xyy} + d\,u_{yyy} +
f(D^{2}u,\,Du,\,u,\,x,\,t)=0,\\
&  & u(x,\,y,\,0)=u_{0}(x,\,y),
\end{eqnarray*}
where $a,$ $b,$ $c,$ $d$ are assumed constant. Indeed, Levandosky
proved sufficient conditions on this equation for which a solution
$u$ will experience an infinite gain in regularity. Specifically,
prove conditions for which initial data $u_{0}(x,\,y)$ possessing
sufficient decay at infinity and a minimal amount of regularity will
lead to a unique solution $u(t)\in C^{\infty}(\mathbb{R}^{2})$ for
$T^{*}$ where $T^{*}$ is the existence time of solutions. According
to the characteristics of equations $(P)$ and considering the
particular cases $(Q_{1})$ and $(Q_{2})$ we could hope that the
$(P)$ equation have gain in regularity following the steps of
N. Hayashi {\it et al.} \cite{ha1} or W. Craig {\it et al.} \cite{cr1}.\\
In our problem, the initial idea is to apply the technique given by
N. Hayashi {\it et al.} \cite{ha1,ha2} to obtain gain in regularity.
Firstly, using straightforward calculus we can see that the equation
$(P)$ has conservation of the energy, i. e.,
$||u||_{L^{2}(\mathbb{R})}=||u_{0}||_{L^{2}(\mathbb{R})}.$ On the
other hand, we look for estimates for $u_{x}$ that will help to
obtain a priori estimates, basically to obtain estimates in
$L^{\infty}(\mathbb{R}).$ Indeed, differentiating in the
$x$-variable the equation $(P)$ we have
\begin{eqnarray}
\label{e100}i\,u_{x\,t} + i\,\beta \,u_{xxxx} + \omega \,u_{xxx} +
(|u|^{2})_{x}\,u + |u|^{2}\,u_{x}=0,
\end{eqnarray}
and multiplying \eqref{e100} by $\overline{u}_{x}$
\begin{eqnarray*}
&  & i\,\overline{u}_{x}\,u_{x\,t} + i\,\beta
\,\overline{u}_{x}\,u_{xxxx} + \omega\,\overline{u}_{x}\,u_{xxx} +
\,(|u|^{2})_{x}\,u\,\overline{u}_{x} + |u|^{2}\,|u_{x}|^{2}=0\\
&  & -\,i\,u_{x}\,\overline{u}_{x\,t} - i\,\beta
\,u_{x}\,\overline{u}_{xxxx} + \omega \,u_{x}\,\overline{u}_{xxx}
+ \,(|u|^{2})_{x}\,\overline{u}\,u_{x} + |u|^{2}\,|u_{x}|^{2}=0.\;
(\mbox{applying conjugate})
\end{eqnarray*}
Subtracting and integrating over $x\in \mathbb{R},$ we have
\begin{eqnarray*}
\lefteqn{i\,\frac{d}{dt}\int_{\mathbb{R}}|u_{x}|^{2}dx + i\,\beta
\int_{\mathbb{R}}\overline{u}_{x}\,u_{xxxx}dx + i\,\beta
\int_{\mathbb{R}}u_{x}\,\overline{u}_{xxxx}dx }\\
&  & +\;2\,i\,\omega \,Im
\int_{\mathbb{R}}\overline{u}_{x}\,u_{xxx}dx + 2\,i\,Im
\int_{\mathbb{R}}(|u|^{2})_{x}\,u\,\overline{u}_{x}dx=0.
\end{eqnarray*}
Performing integration by parts and straightforward calculations we
obtain
\begin{eqnarray*}
\frac{d}{dt}\int_{\mathbb{R}}|u_{x}|^{2}dx + 2\,Im
\int_{\mathbb{R}}(|u|^{2})_{x}\,u\,\overline{u}_{x}dx=0\qquad
(E_{1})
\end{eqnarray*}
where
\begin{eqnarray*}
\frac{d}{dt}\,||u_{x}||_{L^{2}(\mathbb{R})}^{2} + 2\,Im
\int_{\mathbb{R}}u^{2}\,\overline{u}_{x}^{2}dx=0\qquad (E_{2})
\end{eqnarray*}
or integrating by parts the second term in $(E_{1})$ we obtain
\begin{eqnarray*}
\frac{d}{dt}\,||u_{x}||_{L^{2}(\mathbb{R})}^{2} - 2\,Im
\int_{\mathbb{R}}|u|^{2}\,u\,\overline{u}_{xx}dx=0.\qquad (E_{3})
\end{eqnarray*}
Thus it is not possible to estimate in $H^{1}(\mathbb{R}),$
because it appears a second term with two derivatives. The reason
of having an estimate in the derivative is related to Sobolev
embedding. In one spatial dimension we have the embedding
$H^{1}(\mathbb{R})\hookrightarrow L^{\infty}(\mathbb{R}).$ It
seems that the term $i\,\beta\,u_{xxx}$ is crucial. It makes the
two "top" terms look like KdV equation; that is, $u_{t} + u_{xxx}
+ \ldots .$ Of course, the solution is complex, so that the
equation
is like two coupled real KdV equations. \\
This was our motivation to obtain gain in regularity using the idea
of W. Craig {\it et al.} \cite{cr1}. We prove conditions on $(P)$
for which initial data $u_{0}(x)$ possessing sufficient decay at
infinity and a minimal amount of regularity will lead to a unique
solution $u(t)\in C^{\infty}(\mathbb{R})$ for $t>0.$ We use a
technique of nonlinear multipliers, generalizing Kato's original
method, together with ideas of Craig and Goodman \cite{cr0} All the
physically significant dispersive equations and systems known to us
have linear parts displaying this local smoothing property. To
mention only a few, the KdV, Benjamin-Ono, intermediate long wave,
various Boussinesq, and Schr\"{o}dinger equation are included. This
paper is organized as follows: Section 2 outlines briefly the
notation and terminology to be used subsequently. In section 3 we
prove the main inequality. In section 4 we prove an important a
priori estimate. In section 5 we prove a basic-local-in-time
existence and uniqueness theorem. In section 6 we prove a basic
global existence theorem. In section 7 we develop a series of
estimates for solutions of equations $(P)$ in weighted Sobolev
norms. These provide a starting point for the a priori gain of
regularity. In section 8 we prove the
following theorem:\\
{\bf Theorem 1.1}(Main Theorem). {\it Let $|\omega |<3\,\beta ,$
$T>0$ and $u(x,\,t)$ be a solution of $(P)$ in the region
$\mathbb{R} \times [0,\,T]$ such that}
\begin{eqnarray}
\label{e101}u\in L^{\infty}([0,\,T]:\,H^{3}(W_{0\;L\;0}))
\end{eqnarray}
{\it for some $L\geq 2.$ Then }
\begin{eqnarray}
\label{e102}u\in L^{\infty}([0,\,T]:\,H^{3 + l}(W_{\sigma,\,L -
l,\,l}))\cap L^{2}([0,\,T]:\,H^{4 + l}(W_{\sigma,\,L - l - 1,\,l}))
\end{eqnarray}
{\it for all $0\leq l\leq L - 1$ and all $\sigma >0.$}\\
\\
{\it Remark.} We consider the Gauge transformation
\begin{eqnarray}
\label{e104}u(x,\,t) & = & e^{i\,d_{2}\,x + i\,d_{3}\,t}\,v\left(x -
d_{1}\,t,\,t\right)\equiv e^{\theta}\,v\left(\eta,\,\xi\right)
\end{eqnarray}
where $\;\theta = i\,d_{2}\,x + i\,d_{3}\,t,\;$ $\eta=x - d_{1}\,t$
and $\xi =t.$ Then
\begin{eqnarray*}
&  & u_{t} = i\,d_{3}\,e^{\theta}\,v - d_{1}\,e^{\theta}\,v_{\eta} +
e^{\theta}\,v_{\xi}\quad :\quad
u_{x} = i\,d_{2}\,e^{\theta}\,v + e^{\theta}\,v_{\eta}\\
&  & u_{xx} = -\;d_{2}^{2}\,e^{\theta}\,v + 2\,i\,d_{2}
\,e^{\theta}\,v_{\eta} + e^{\theta}\,v_{\eta \,\eta}\; :\; u_{xxx} =
-\;i\,d_{2}^{3}\,e^{\theta}\,v - 3\,d_{2}^{2}\,e^{\theta}\,v_{\eta}
+ 3\,i\,d_{2} \,e^{\theta}\,v_{\eta\eta} +
e^{\theta}\,v_{\eta\eta\eta}.
\end{eqnarray*}
Replacing in $(Q)$ we have
\begin{eqnarray*}
&  & -\,d_{3}\,e^{\theta}\,v - i\,d_{1}\,e^{\theta}\,v_{\eta} +
i\,e^{\theta}\,v_{\xi} - \omega\,d_{2}^{2}\,e^{\theta}\,v +
2\,i\,\omega\,d_{2}\,e^{\theta}\,v_{\eta} +
\omega\,e^{\theta}\,v_{\eta\eta} \\
&  & \beta\,d_{3}^{3}\,e^{\theta}\,v -
3\,i\,\beta\,d_{2}^{2}\,e^{\theta}\,v_{\eta} -
3\,\beta\,d_{2}\,e^{\theta}\,v_{\eta\eta} +
i\,\beta\,e^{\theta}\,v_{\eta\eta\eta} + \gamma\,|v|^{2}\,e^{\theta}\,v\\
&  & -\;\delta\,d_{2}\,|v|^{2}\,e^{\theta}\,v +
i\,\delta\,|v|^{2}\,e^{\theta}\,v_{\eta} +
\epsilon\,d_{2}\,e^{\theta}\,v^{2}\overline{v} +
i\,\epsilon\,e^{\theta}\,v^{2}\,v_{\eta}=0
\end{eqnarray*}
where
\begin{eqnarray*}
&  & i\,v_{\xi} + (\omega - 3\,\beta\,d_{2})\,v_{\eta\eta} +
i\,\beta\,v_{\eta\eta\eta} + (2\,i\,\omega\,d_{2} -
3\,i\,\beta\,d_{2}^{2} - i\,d_{1} + i\,\delta\,|v|^{2} +
i\,\epsilon\,v^{2})\,v_{\eta}\\
&  & (\beta\,d_{2}^{3} - \omega\,d_{2}^{2} - d_{3} + \gamma\,|v|^{2}
- \delta\,d_{2}\,|v|^{2})\,v + \epsilon\,d_{2}\,v^{2}\overline{v}=0
\end{eqnarray*}
then
\begin{eqnarray}
\label{e105}d_{1}=\frac{\omega^{2}}{3\,\beta}\quad :\quad
d_{2}=\frac{\omega}{3\,\beta}\quad :\quad
d_{3}=\frac{-\,2\omega^{3}}{27\,\beta^{2}}.
\end{eqnarray}
This way in $(Q)$ we obtain
\begin{eqnarray*}
i\,v_{\xi} + i\,\beta\,v_{\eta\eta\eta} + i\,(\delta\,|v|^{2} +
\epsilon\,v^{2})\,v_{\eta} + \left(\gamma -
\frac{\omega\,\delta}{3\,\beta}\right)|v|^{2}v +
\frac{\epsilon\,\delta}{3\,\beta}\,v^{2}\overline{v}=0,
\end{eqnarray*}
but $v^{2}\,\overline{v}=v\,v\,\overline{v}=|v|^{2}v,$ then using
the Gauge transformation we have the equivalent problem to $(Q)$
\[(\mathbb{Q}) \left \{ \begin{array}{ll}
i\,v_{\xi} + i\,\beta\,v_{\eta\eta\eta} +
i\,\delta\,|v|^{2}\,v_{\eta} + i\,\epsilon\,v^{2}\,v_{\eta} +
\left(\gamma + \frac{\epsilon\,\delta}{3\,\beta} -
\frac{\omega\,\delta}{3\,\beta}\right)|v|^{2}v =0 &
\quad \eta,\,\xi\in \mathbb{R}\\
v(\eta,\,0)  = e^{-\,i\,\frac{\omega}{3\,\beta}\,\eta}u_{0}(\eta). &
\end{array}
\right. \] Here, rescaling the equation, we take $\beta =1.$
\[(\widetilde{\mathbb{Q}}) \left \{ \begin{array}{ll}
i\,v_{t} + i\,v_{xxx} + i\,\delta\,|v|^{2}\,v_{x} +
i\,\epsilon\,v^{2}\,v_{x} + \left(\gamma +
\frac{\epsilon\,\delta}{3} - \frac{\omega\,\delta}{3}\right)|v|^{2}v
=0 &
\quad x,\,t\in \mathbb{R}\\
v(x,\,0)  = e^{-\,i\,\frac{\omega}{3}\,x}u_{0}(x). &
\end{array}
\right. \] The above Gauge transformation is a bicontinuous map
from $L^{p}([0,\,T]:\,H^{s}(W_{\sigma\,i\,k}))$ to itself, as far
as $0<T<+\infty$ and $p,$ $s,$ $\sigma,$ $i,$ $k$ used in this
paper. With this, the assumption $|\omega|<3\,\beta$ imposed in
Theorem 1.1 can be removed.
\renewcommand{\theequation}{\thesection.\arabic{equation}}
\setcounter{equation}{0}\section{Preliminaries}
We consider the initial value problem
\[(P)\left \{ \begin{array}{ll}
i\,u_{t} + \omega \,u_{xx} + i\,\beta \,u_{xxx} + |u|^{2}\,u=0, &
\quad x,\,t\in \mathbb{R} \\
u(x,\,0)  = u_{0}(x) &
\end{array}
\right. \] where $\omega ,\,\beta \in \mathbb{R},$ $\beta \neq 0 $
and $u=u(x,\,t)$ is a
complex valued function.\\
\\
{\it Notation.} We write $\;\partial =\partial /\partial x,\;$
$\;\partial _{t}=\partial /\partial t\;$ and we abbreviate
$\;u_{j}=\partial ^{j}u.$\\
\\
{\it Definition 2.1.} A function $\xi = \xi (x,\,t)$ belongs to
the weight class $W_{\sigma \;i\; k}$ if it is a positive
$C^{\infty }$ function on $\mathbb{R}\times [0,\,T],$ $\partial
\xi
>0$ and there are constant $c_{j},$ $0\leq j\leq 5$ such that
\begin{eqnarray}
\label{e201}&  & 0<c_{1}\leq t^{-\,k}\,e^{-\,\sigma \,x}\,\xi
(x,\,t) \leq c_{2 }\qquad \forall
\;x<-1,\quad 0<t<T.\\
\label{e202}&  & 0<c_{3}\leq t^{-\,k}\,x^{-\,i}\,\xi (x,\,t)\leq
c_{4}\qquad \forall
\;x>1,\quad 0<t<T.\\
\label{e203}&  & \left(t\mid \partial_{t}\xi \mid + \mid
\partial^{j}\xi \mid \right)/ \xi \leq c_{5}\quad \forall
\;(x,\,t)\in \mathbb{R}\times [0,\,T],\;\forall\;j\in\mathbb{N}.
\end{eqnarray}
{\it Remark.}  We shall always take $\sigma \geq 0,$ $i\geq 1$
and $k\geq 0.$\\
\\
{\it Example.} Let
\[ \xi (x)= \left \{\begin{array}{ll}
1 + e^{-1/x} & \mbox { for $\,x>0$ }   \\
1 & \mbox { for $x\leq 0$ }
\end{array}
\right. \]
then $\xi \in W_{0\;i\;0}.$\\
\\
{\it Definition 2.2.} Let $N$ be a positive integer. By
$H^{N}(W_{\sigma \;i\; k})$ we denote the Sobolev space on
$\mathbb{R}$ with a weight; that is, with the norm
\begin{eqnarray*}
||v||_{H^{N}(W_{\sigma \;i\;k})}^{2}= \sum _{j=0}^{N}
\int_{\mathbb{R}}|
\partial ^{j}v(x)|^{2}\,\xi (x,\,t)\,dx<+\,\infty
\end{eqnarray*}
for any $\xi \in W_{\sigma \;i\;k}$ and $0<t<T.,$
Even though the norm depends on $\xi,$ all such choices leads to equivalent norms. \\
\\
{\it Remark.} $H^{s}(W_{\sigma \;i\;k})\,$ depends on $t$ (because
$\xi =\xi (x,\,t)$).\\
\\
{\bf Lemma 2.1.} (See \cite{ca0}) For $\xi \in W_{\sigma \;i\;0}$
and $\sigma \geq 0,\, i\geq 0,$ there exists a constant $c>0$ such
that, for $u\in H^{1}(W_{\sigma \;i\;0}),$
\[\sup _{x\in \mathbb{R}}||\xi \,u^{2}||
\leq c\int _{\mathbb{R}}\left(\,|u|^{2} + |
\partial u|^{2}\,\right)\,\xi \,dx\]
{\bf Lemma 2.2}(The Gagliardo-Nirenberg inequality). Let $q,\,r$
be any real numbers satisfying $1\leq q,$ $r\leq \infty $ and let
$j$ and $m$ be nonnegative integers such that $j\leq m.$ Then
\begin{eqnarray*}
||\partial^{j}u||_{L^{p}(\mathbb{R})}\leq
c\;||\partial^{m}u||_{L^{r}(\mathbb{R})}^{a}\;||u||_{L^{q}(\mathbb{R})}^{1
- a}
\end{eqnarray*}
where $\frac{1}{p}=j + a\,\left (\frac{1}{r} - m\right ) +
\frac{(1 - a)}{q}$ for all $a$ in the interval $\frac{j}{m}\leq a\leq 1,$ and
$M$ is a positive constant depending only on $m,$ $j,$ $q,$ $r$ and $a.$\\
\\
{\it Definition 2.3.} By $L^{2}([0,\,T]:\,H^{N}(W_{\sigma \;i\;k}))$
we denote the space of functions $v(x,\,t)$ with the norm ($N$
integer positive)
\begin{eqnarray*}
||v||_{L^{2}([0,\,T]:\,H^{N}(W_{\sigma \;i\;k}))}^{2}=
\int_{0}^{T}||v(x,\,t)||_{H^{N}(W_{\sigma
\;i\;k})}^{2}dt<+\,\infty
\end{eqnarray*}
{\it Remark.} The usual Sobolev space is $\,H^{N}(\mathbb{R}) =
H^{N}(W_{0\;0\;0})\,$
without a weight.\\
\\
{\it Remark.} We shall derive the a priori estimates assuming that
the solution is $C^{\infty},$ bounded as $x\rightarrow -\,\infty ,$
and rapidly decreasing as
$x\rightarrow +\,\infty ,$ together with all of its derivatives.\\
\\
Considering the above notation, the higher order nonlinear
Schr\"{o}dinger equation can be written as
\begin{eqnarray}
\label{e204}i\,u_{t} + i\,\beta\,u_{3} + \omega \,u_{2} +
|u|^{2}\,u=0,\quad x,\,t\in \mathbb{R}
\end{eqnarray}
where $\omega ,\,\beta \in \mathbb{R},$ $\beta \neq 0 $ and
$u=u(x,\,t)$ is a
complex valued function.\\
\\
Throughout this paper  $c$ is a generic constant, not necessarily
the same at each occasion(it will change from line to line), which
depends in an increasing way on the indicated quantities. In this
part, we only consider the case $t>0.$ The case $t<0$ can be
treated analogously.
\renewcommand{\theequation}{\thesection.\arabic{equation}}
\setcounter{equation}{0}\section{Main Inequality} {\bf Lemma 3.1.}
{\it Let $|\omega |<3\;\beta .$ Let $u$ be a solution of
\eqref{e204} with enough Sobolev regularity (for instance, $u\in
H^{N}(\mathbb{R}),$ $N\geq \alpha + 3$), then}
\begin{eqnarray}
\label{e301}&  & \partial_{t}\int_{\mathbb{R}}\xi
\,|u_{\alpha}|^{2}dx + \int_{\mathbb{R}}\eta \,|u_{\alpha +
1}|^{2}dx + \int_{\mathbb{R}}\theta \,|u_{\alpha}|^{2}dx +
\int_{\mathbb{R}}R_{\alpha}dx\leq 0
\end{eqnarray}
{\it where}
\begin{eqnarray*}
\eta & = & (3\,\beta - |\omega |)\,\partial \xi\qquad for\qquad
|\omega |<3\;
\beta \\
\theta & = & -\;[\,\partial_{t}\xi + \beta \,\partial^{3}\xi +
|\omega|\,\partial \xi + c_{0}\,\xi \,] \quad where\quad
c_{0}=||u||_{L^{\infty}(\mathbb{R})}^{2}
\end{eqnarray*}
{\it and $R_{\alpha}=R_{\alpha}(|u_{\alpha }|,\,|u_{\alpha - 1}|,\,\ldots ).$}\\
\\
{\it Proof.} Differentiating \eqref{e204} $\alpha $-times (for
$\alpha \geq 0$) over $x\in \mathbb{R}$ leads to
\begin{eqnarray}
\label{e302}i\,u_{\alpha \,t} + i\,\beta \,u_{\alpha + 3} + \omega
\,u_{\alpha + 2} + (|u|^{2})_{\alpha }\,u + \sum_{m=1}^{\alpha -
1}{\alpha\choose m}\,(|u|^{2})_{\alpha - m}\,u_{m} +
|u|^{2}\,u_{\alpha } = 0.
\end{eqnarray}
Let $\xi = \xi(x,\,t),$ then multiplying \eqref{e302} by $\xi
\,\overline{u}_{\alpha } $
 we have
\begin{eqnarray*}
\lefteqn{i\,\xi \,\overline{u}_{\alpha }\,u_{\alpha \,t} +
i\,\beta \,\xi \, \overline{u}_{\alpha }\,u_{\alpha + 3} + \omega
\,\xi \,\overline{u}_{\alpha }\,
u_{\alpha + 2} + (|u|^{2})_{\alpha }\,\xi \,u\,\overline{u}_{\alpha } } \\
&  & +\sum_{m=1}^{\alpha - 1}{\alpha\choose m}\,(|u|^{2})_{\alpha -
m}\,
\xi \,u_{m}\,\overline{u}_{\alpha} + \xi \,|u|^{2}\,|u_{\alpha }|^{2}  =  0 \\
\\
\lefteqn{-\,i\,\xi \,u_{\alpha }\,\overline{u}_{\alpha \,t} -
i\,\beta \,\xi \,u_{\alpha }\, \overline{u}_{\alpha + 3} + \omega
\,\xi \,u_{\alpha }\,\overline{u}_{\alpha + 2} +
(|u|^{2})_{\alpha }\,\xi \,\overline{u}\,u_{\alpha } } \\
&  & +\sum_{m=1}^{\alpha - 1}{\alpha\choose m}\,(|u|^{2})_{\alpha -
m} \,\xi \,\overline{u}_{m}\,u_{\alpha} + \xi \,|u|^{2}\,|u_{\alpha
}|^{2}  =  0 . \qquad (\mbox{applying conjugate})
\end{eqnarray*}
Subtracting and integrating over $x\in \mathbb{R}$ we have
\begin{eqnarray}
\lefteqn{i\,\partial_{t}\int_{\mathbb{R}}\xi \,|u_{\alpha }|^{2}dx
+ i\,\beta \int_{\mathbb{R}}\xi \,\overline{u}_{\alpha
}\,u_{\alpha + 3}dx + i\,\beta \int_{\mathbb{R}}\xi
\,u_{\alpha}\,\overline{u}_{\alpha + 3}dx
- i\int_{\mathbb{R}}\xi_{t}\,|u_{\alpha }|^{2}dx } \nonumber \\
&  & +\;\omega \int_{\mathbb{R}}\xi
\,\overline{u}_{\alpha}\,u_{\alpha + 2}dx - \omega
\int_{\mathbb{R}}\xi \,u_{\alpha}\,\overline{u}_{\alpha + 2}dx +
2\,i\,Im\int_{\mathbb{R}}\xi \,(|u|^{2})_{\alpha
}\,u\,\overline{u}_{\alpha }dx \nonumber \\
\label{e303}&  & +\;2\,i\sum_{m=1}^{\alpha - 1}{\alpha\choose
m}\;Im\int_{\mathbb{R}}\xi \,(|u|^{2})_{\alpha - m}
\,u_{m}\,\overline{u}_{\alpha}dx = 0.
\end{eqnarray}
We estimate the second term integrating by parts
\begin{eqnarray*}
&  & \int_{\mathbb{R}}\xi\,\overline{u}_{\alpha}\,u_{\alpha + 3}dx =
\int_{\mathbb{R}}\partial^{2}\xi \,\overline{u}_{\alpha}\,u_{\alpha
+ 1}dx + 2\int_{\mathbb{R}}\partial\xi \,|u_{\alpha + 1}|^{2}dx +
\int_{\mathbb{R}}\xi\,\overline{u}_{\alpha + 2}\,u_{\alpha + 1}dx.
\end{eqnarray*}
The other terms are calculated in a similar way. Hence, replacing in
\eqref{e303}  and performing straightforward calculations we obtain
\begin{eqnarray*}
\lefteqn{i\,\partial_{t}\int_{\mathbb{R}}\xi \,|u_{\alpha }|^{2}dx +
i\,\beta
\int_{\mathbb{R}}\partial^{2}\xi\,\overline{u}_{\alpha}\,u_{\alpha +
1}dx + 2\,i\,\beta \int_{\mathbb{R}}\partial \xi\,|u_{\alpha +
1}|^{2}dx }\\
&  & +\;i\,\beta \int_{\mathbb{R}}\xi \,\overline{u}_{\alpha +
2}\,u_{\alpha + 1}dx + i\,\beta \int_{\mathbb{R}}\partial^{2}\xi
\,u_{\alpha}\,\overline{u}_{\alpha + 1}dx
 +  i\,\beta \int_{\mathbb{R}}\partial \xi\,|u_{\alpha + 1}|^{2}dx \\
 &  & -\;
 i\,\beta \int_{\mathbb{R}}\xi\,u_{\alpha + 1}\,\overline{u}_{\alpha + 2}dx -
 \omega \int_{\mathbb{R}}\partial \xi\,\overline{u}_{\alpha}\,u_{\alpha + 1}dx -
 \omega \int_{\mathbb{R}}\xi \,|u_{\alpha + 1}|^{2}dx \\
 &  & +\;
 \omega \int_{\mathbb{R}}\partial \xi\,u_{\alpha}\,\overline{u}_{\alpha + 1}dx
 +  \omega \int_{\mathbb{R}}\xi \,|u_{\alpha + 1}|^{2}dx -
 i\int_{\mathbb{R}}\partial_{t}\xi\,|u_{\alpha}|^{2}dx \\
 &  & +\;2\,i\,Im\int_{\mathbb{R}}\xi \,(|u|^{2})_{\alpha }\,
 u\,\overline{u}_{\alpha }dx +
 2\,i\sum_{m=1}^{\alpha - 1}{\alpha\choose m}Im\int_{\mathbb{R}}\xi \,
 (|u|^{2})_{\alpha - m} \,u_{m}\,\overline{u}_{\alpha}dx = 0
\end{eqnarray*}
then
\begin{eqnarray*}
\lefteqn{\partial_{t}\int_{\mathbb{R}}\xi \,|u_{\alpha}|^{2}dx -
\beta \int_{\mathbb{R}}\partial^{3}\xi\,|u_{\alpha}|^{2}dx +
3\,\beta \int_{\mathbb{R}}\partial \xi\,|u_{\alpha + 1}|^{2}dx -
2\,\omega \,Im\int_{\mathbb{R}}\partial \xi\,
\overline{u}_{\alpha}\,u_{\alpha + 1}dx}\\
&  & -\int_{\mathbb{R}}\partial_{t}\xi\,|u_{\alpha}|^{2}dx +
2\,Im\int_{\mathbb{R}}\xi\,(|u|^{2})_{\alpha
}\,u\,\overline{u}_{\alpha}dx + 2\sum_{m=1}^{\alpha -
1}{\alpha\choose m}Im\int_{\mathbb{R}}\xi \,(|u|^{2})_{\alpha - m}
\,u_{m}\,\overline{u}_{\alpha}dx = 0
\end{eqnarray*}
hence
\begin{eqnarray*}
\lefteqn{\partial_{t}\int_{\mathbb{R}}\xi\,|u_{\alpha}|^{2}dx -
\beta \int_{\mathbb{R}}\partial^{3}\xi\,|u_{\alpha}|^{2}dx +
3\,\beta \int_{\mathbb{R}}\partial \xi\,|u_{\alpha + 1}|^{2}dx +
2\,Im\int_{\mathbb{R}}(|u|^{2})_{\alpha}\,\xi
\,u\,\overline{u}_{\alpha}dx}\\
&  & -\int_{\mathbb{R}}\partial_{t}\xi\,|u_{\alpha}|^{2}dx
+\,2\sum_{m=1}^{\alpha - 1}{\alpha\choose
m}Im\int_{\mathbb{R}}\xi\,(|u|^{2})_{\alpha - m}
\,u_{m}\,\overline{u}_{\alpha}dx =
2\,\omega\,Im\int_{\mathbb{R}}\partial \xi\,
\overline{u}_{\alpha}\,u_{\alpha + 1}dx \\
& \leq & |\omega |\int_{\mathbb{R}}\partial \xi\,|u_{\alpha}|^{2}dx
+ |\omega |\int_{\mathbb{R}}\partial \xi\,|u_{\alpha + 1}|^{2}dx
\end{eqnarray*}
therefore
\begin{eqnarray}
\lefteqn{\partial_{t}\int_{\mathbb{R}}\xi \,|u_{\alpha}|^{2}dx +
\int_{\mathbb{R}}[\,3\,\beta - |\omega |\,]\,\partial
\xi\;|u_{\alpha + 1}|^{2}dx - \int_{\mathbb{R}}[\,\partial_{t}\xi +
\beta \,\partial^{3}\xi + |\omega |\,\partial
\xi\,]\,|u_{\alpha}|^{2}dx }\nonumber \\
\label{e304}&  & +\;2\,Im\int_{\mathbb{R}}(|u|^{2})_{\alpha}\,\xi
\,u\,\overline{u}_{\alpha}dx + 2\sum_{m=1}^{\alpha -
1}{\alpha\choose m}Im\int_{\mathbb{R}}\xi \,(|u|^{2})_{\alpha - m}
\,u_{m}\,\overline{u}_{\alpha}dx\leq 0.
\end{eqnarray}
But
\begin{eqnarray*}
(|u|^{2})_{\alpha} & = & (\,u\,\overline{u}\,)_{\alpha}  =
\sum_{k=0}^{\alpha}{\alpha\choose k}u_{\alpha - k}\,
\overline{u}_{k} = \overline{u}\,u_{\alpha} + \sum_{k=1}^{\alpha -
1}{\alpha\choose k}u_{\alpha - k}\,\overline{u}_{k} +
u\,\overline{u}_{\alpha}
\end{eqnarray*}
then
\begin{eqnarray*}
(|u|^{2})_{\alpha}\,u\,\overline{u}_{\alpha} =
|u|^{2}|u_{\alpha}|^{2} + \sum_{k=1}^{\alpha - 1}{\alpha\choose
k}u_{\alpha - k}\, \overline{u}_{k}\,u\,\overline{u}_{\alpha} +
u^{2}\,\overline{u}_{\alpha}^{2}
\end{eqnarray*}
thus,
\begin{eqnarray}
\lefteqn{2\;Im\int_{\mathbb{R}}(\,|u|^{2})_{\alpha}\,\xi
\,u\,\overline{u}_{\alpha}dx = 2\sum_{k=1}^{\alpha -
1}{\alpha\choose k}Im\int_{\mathbb{R}}\xi\,u_{\alpha -
k}\,\overline{u}_{k}\,u\,\overline{u}_{\alpha}dx +
2\;Im\int_{\mathbb{R}}\xi \,
u^{2}\,\overline{u}_{\alpha}^{2}dx }\nonumber \\
& \leq & 2\sum_{k=1}^{\alpha - 1}{\alpha\choose
k}\int_{\mathbb{R}}\xi\,|u_{\alpha -
k}|\,|u_{k}|\,|u|\,|u_{\alpha}|dx + 2\int_{\mathbb{R}}\xi
\,|u|^{2}\,
|u_{\alpha}|^{2}dx \nonumber \\
& \leq & 2\sum_{k=1}^{\alpha - 1}{\alpha\choose
k}\int_{\mathbb{R}}\xi\,|u_{\alpha -
k}|\,|u_{k}|\,|u|\,|u_{\alpha}|dx +
2\,||u||_{L^{\infty}(\mathbb{R})}^{2}\int_{\mathbb{R}}\xi
\,|u_{\alpha}|^{2}dx \nonumber \\
\label{e305}& \leq &
2\,||u||_{L^{\infty}(\mathbb{R})}\sum_{k=1}^{\alpha -
1}{\alpha\choose k}\int_{\mathbb{R}}\xi \,|u_{\alpha -
k}|\,|u_{k}|\,|u_{\alpha}|dx +
2\,||u||_{L^{\infty}(\mathbb{R})}^{2}\int_{\mathbb{R}}\xi
\,|u_{\alpha}|^{2}dx\qquad \;
\end{eqnarray}
hence, in \eqref{e304} we have
\begin{eqnarray*}
\lefteqn{\partial_{t}\int_{\mathbb{R}}\xi \,|u_{\alpha}|^{2}dx +
\int_{\mathbb{R}}[3\,\beta - |\omega |\,]\,\partial \xi\,|u_{\alpha
+ 1}|^{2}dx - \int_{\mathbb{R}}[\partial_{t}\xi + \beta
\,\partial^{3}\xi +
|\omega|\,\partial \xi + c_{0}\,\xi \,]\, |u_{\alpha}|^{2}dx} \\
&  & -\;2\,c\sum_{k=1}^{\alpha - 1}{\alpha\choose
k}\int_{\mathbb{R}}\xi \,|u_{\alpha - k}|\,|u_{k}|\,|u_{\alpha}|dx -
2\sum_{m=1}^{\alpha - 1}{\alpha\choose m}\int_{\mathbb{R}}\xi \,
|(|u|^{2})_{\alpha - m}| \,|u_{m}|\,|u_{\alpha}|dx\leq 0.
\end{eqnarray*}
Therefore, using straightforward calculations we obtain the {\it
main inequality}
\begin{eqnarray*}
&  & \partial_{t}\int_{\mathbb{R}}\xi \,|u_{\alpha}|^{2}dx +
\int_{\mathbb{R}}\eta \,|u_{\alpha + 1}|^{2}dx +
\int_{\mathbb{R}}\theta \,|u_{\alpha}|^{2}dx +
\int_{\mathbb{R}}R_{\alpha}dx\leq 0
\end{eqnarray*}
where
\begin{eqnarray*}
\eta & = & (3\,\beta - |\omega |\,)\,\partial \xi\qquad
\mbox{for}\qquad
|\omega |<3\;\beta \\
\theta & = & -\;[\,\partial_{t}\xi + \beta \,\partial^{3}\xi +
|\omega|\,\partial \xi + c_{0}\,\xi \,] \quad \mbox{where}\quad
c_{0}=||u||_{L^{\infty}(\mathbb{R})}^{2}
\end{eqnarray*}
and $R_{\alpha}=R_{\alpha}(|u_{\alpha }|,\,|u_{\alpha - 1}|,\,\ldots ).$\\
\\
{\it Remark.} In \eqref{e304} using Young's estimate and assuming
that $\beta>0$ we have
\begin{eqnarray*}
2\,\omega\;Im\int_{\mathbb{R}}\overline{u_{\alpha}}\,u_{\alpha +
1}\,dx\leq
\frac{|\omega|^{2}}{2\,\beta}\int_{\mathbb{R}}|u_{\alpha}|^{2}\,dx +
2\,\beta\int_{\mathbb{R}}|u_{\alpha + 1}|^{2}\,dx.
\end{eqnarray*}
Then, in \eqref{e304} we obtain
\begin{eqnarray*}
\lefteqn{\partial_{t}\int_{\mathbb{R}}\xi\,|u_{\alpha}|^{2}dx -
\beta \int_{\mathbb{R}}\partial^{3}\xi\,|u_{\alpha}|^{2}dx + \beta
\int_{\mathbb{R}}\partial \xi\,|u_{\alpha + 1}|^{2}dx +
2\,Im\int_{\mathbb{R}}(|u|^{2})_{\alpha}\,\xi
\,u\,\overline{u}_{\alpha}dx}\\
&  & -\int_{\mathbb{R}}\partial_{t}\xi\,|u_{\alpha}|^{2}dx
+\,2\sum_{m=1}^{\alpha - 1}{\alpha\choose
m}Im\int_{\mathbb{R}}\xi\,(|u|^{2})_{\alpha - m}
\,u_{m}\,\overline{u}_{\alpha}dx =
2\,\omega\,Im\int_{\mathbb{R}}\partial \xi\,
\overline{u}_{\alpha}\,u_{\alpha + 1}dx \\
& \leq &
\frac{|\omega|^{2}}{2\,\beta}\int_{\mathbb{R}}|u_{\alpha}|^{2}\,dx
\end{eqnarray*}
and the assumption that $|\omega|<3\,\beta$ can be removed.\\
\\
{\bf Lemma 3.2.} {\it For $\eta \in W_{\sigma \;i\;k}$ an arbitrary
weight function and $|\omega |<3\,\beta ,$ there exists $\xi \in
W_{\sigma ,\;i + 1,\;k}$ that satisfies}
\begin{eqnarray}
\label{e306}\eta = (3\,\beta - |\omega |)\,\partial \xi \qquad
for\qquad |\omega |<3\;\beta .
\end{eqnarray}
Indeed, we have
\begin{eqnarray}
\label{e307}\xi = \frac{1}{(3\,\beta - |\omega
|)}\int_{-\,\infty}^{x}\eta (y,\,t)\,dy.
\end{eqnarray}
{\bf Lemma 3.3.} {\it The expression $R_{\alpha}$ in the inequality
of Lemma 3.1 is a sum of terms of the form}
\begin{eqnarray}
\label{e308}\xi
\,u_{\nu_{1}}\;\overline{u}_{\nu_{2}}\;\overline{u}_{\alpha}
\end{eqnarray}
{\it where $1\leq \nu_{1}\leq \nu_{2}\leq \alpha$ and }
\begin{eqnarray}
\label{e309}\nu_{1} + \nu_{2} = \alpha
\end{eqnarray}
{\it Proof.} It follows from \eqref{e305}.
\renewcommand{\theequation}{\thesection.\arabic{equation}}
\setcounter{equation}{0}\section{An a priori estimate} We show now a
fundamental a priori estimate used for a basic local-in-time
existence theorem. We construct a mapping ${\cal
Z}:L^{\infty}([0,\,T]:\,H^{s}(\mathbb{R}))\longmapsto
L^{\infty}([0,\,T]:\,H^{s}(\mathbb{R}))$ with the property:\\
Given $u^{(n)}={\cal Z}(u^{(n - 1)})$ and
$essup_{t\in[0,\,T]}||u^{(n - 1)}||_{s}\leq c_{0}$ then
$essup_{t\in[0,\,T]}||u^{(n)}||_{s}\leq c_{0},$ where $s$ and
$c_{0}>0$ are constants. This property tells us that ${\cal
Z}:\mathbb{B}_{c_{0}}(0)\longmapsto \mathbb{B}_{c_{0}}(0)$ where
$\mathbb{B}_{c_{0}}(0)=\{v(x,\,t):\;||v(\,\cdot\,,\,t)||_{s}\leq
c_{0}\}$ is a ball in $L^{\infty}([0,\,T]:\,H^{s}(\mathbb{R})).$
To guarantee this property, we will appeal to an a
priori estimate which is the main object of this section.\\
\\
Differentiating \eqref{e204} two times leads to
\begin{eqnarray}
\label{e401}i\,\partial_{t} u_{2} + i\,\beta \,u_{5} + \omega
\,u_{4} + (|u|^{2})_{2}\,u + 2\,(|u|^{2})_{1}\,u_{1} +
|u|^{2}\,u_{2} = 0.
\end{eqnarray}
Let $u=\wedge v$ where $\wedge =(I - \partial^{2})^{-1}.$ Hence
$u=(I - \partial^{2})^{-1}v\;$ then $\;u - u_{2} = v\;$ where
$\;\partial _{t}u_{2} = -\,v_{t} + u_{t}.$ \\
\\
Replacing in \eqref{e401} we have
\begin{eqnarray}
\lefteqn{-\,i\,v_{t} + i\,\beta \,\wedge v_{5} + \omega \,\wedge
v_{4} + (|\wedge v|^{2})_{2}\wedge v + 2\,(|\wedge v|^{2})_{1}\wedge
v_{1} }
\nonumber \\
\label{e402}&  & +\;|\wedge v|^{2}\,\wedge v_{2} - (i\,\beta
\,\wedge v_{3} + \omega \,\wedge v_{2} + |\wedge v|^{2}\,\wedge
v)=0.
\end{eqnarray}
The \eqref{e402} equation is linearized by substituting a new
variable $z$ in each coefficient:
\begin{eqnarray}
\lefteqn{-\,i\,v_{t} + i\,\beta \,\wedge v_{5} + \omega \,\wedge
v_{4} + (|\wedge z|^{2})_{2}\wedge v + 2\,(|\wedge z|^{2})_{1}\wedge
v_{1} }
\nonumber \\
\label{e403}&  & +\;|\wedge z|^{2}\,\wedge v_{2} - (i\,\beta
\,\wedge v_{3} + \omega \,\wedge v_{2} + |\wedge z|^{2}\,\wedge
v)=0.
\end{eqnarray}
The linear equation which is to be solved at each iteration is of
the form
\begin{eqnarray}
\label{e404}i\,\partial_{t}v=i\,\beta \,\wedge v_{5}^{(n)} + \omega
\,\wedge v_{4}^{(n)} - i\,\beta \,\wedge v_{3}^{(n)} - \omega
\,\wedge v_{2}^{(n)} + b^{(1)}
\end{eqnarray}
where $b^{(1)}=(|\wedge z|^{2})_{2}\,\wedge v + 2\,(|\wedge
z|^{2})_{1}\,\wedge v_{1} + |\wedge z|^{2}\,\wedge v_{2} - |\wedge
z|^{2}\,\wedge v.$ Equation \eqref{e404} is a linear equation at
each iteration which can be solved in any interval of time in which
the
coefficient is defined.\\
\\
We consider the following lemma that will help us setting up the
iteration scheme.\\
\\
{\bf Lemma 4.1.} {\it Let $|\omega |<3\,\beta .$ Given initial data
$u_{0}(x)\in H^{\infty}(\mathbb{R}) = \bigcap _{N\geq
0}H^{N}(\mathbb{R})$ there exists a unique solution of \eqref{e404}
where $b^{(1)}$ is a smooth bounded coefficient with $z\in
H^{\infty}(\mathbb{R}).$ The solution is defined in any time
interval in which
the coefficient is defined.}\\
\\
{\it Proof.} Let $T>0$ be arbitrary and $M>0$ a constant. Let
\begin{eqnarray*}
\Gamma =\xi \,(\,i\,\partial_{t} - i\,\beta \,\wedge \partial^{5} -
\omega \,\wedge \partial^{4} + i\,\beta \,\wedge \partial^{3} +
\omega \,\wedge \partial^{2}\,)
\end{eqnarray*}
then in \eqref{e404} we have $\Gamma u=\xi \,b^{(1)}.$  We consider
the bilinear form ${\cal B} :{\cal D}\times {\cal D}\longmapsto
\mathbb{R},$
\begin{eqnarray*}
{\cal
B}(u,\,v)=<u,\,v>=Im\int_{0}^{T}\int_{\mathbb{R}}e^{-Mt}\,u\,
\overline{v}\,dx\,dt
\end{eqnarray*}
where ${\cal D}=\{u\in C_{0}^{\infty}(\mathbb{R} \times
[0,\,T]):\;u(x,\,0)=0\,\}.$ We have
\begin{eqnarray*}
\Gamma u\cdot \overline{u} & = & i\,\xi \,\overline{u}\,u_{t} -
i\,\beta \,\xi \, \overline{u}\,\wedge u_{5} - \omega \,\xi
\,\overline{u}\,\wedge u_{4} + i\,\beta \,\xi \,
\overline{u}\,\wedge u_{3} + \omega \,\xi \,\overline{u}\,\wedge u_{2}\\
\overline{\Gamma u\cdot \overline{u}} & = & -\,i\,\xi
\,u\,\overline{u}_{t} + i\,\beta \, \xi \,u\,\wedge \overline{u}_{5}
- \omega \,\xi \,u\,\wedge \overline{u}_{4} - i\,\beta \,\xi
\,u\,\wedge \overline{u}_{3} + \omega \,\xi \,u\,\wedge
\overline{u}_{2}.\; (\mbox{applying conjugate})
\end{eqnarray*}
Subtracting and integrating over $x\in \mathbb{R}$ we have
\begin{eqnarray*}
\lefteqn{2\,i\,Im\int_{\mathbb{R}}\Gamma u\cdot \overline{u}dx
 =  i\,\partial_{t}\int_{\mathbb{R}}\xi \,|u|^{2}dx -
i\int_{\mathbb{R}}\partial_{t}\xi\,|u|^{2}dx - i\,\beta
\int_{\mathbb{R}}\xi \,\overline{u}\,\wedge u_{5}dx -
i\,\beta \int_{\mathbb{R}}\xi \,u\,\wedge \overline{u}_{5}dx} \\
&  & -\,\omega \int_{\mathbb{R}}\xi \,\overline{u}\,\wedge u_{4}dx +
\omega \int_{\mathbb{R}}\xi \,u\,\wedge \overline{u}_{4}dx +
i\,\beta \int_{\mathbb{R}}\xi \,\overline{u}\,\wedge u_{3}dx +
i\,\beta \int_{\mathbb{R}}\xi \,u\,\wedge \overline{u}_{3}dx \\
&  & +\,\omega \int_{\mathbb{R}}\xi \,\overline{u}\,\wedge u_{2}dx -
\omega \int_{\mathbb{R}}\xi \,u\,\wedge \overline{u}_{2}dx.
\end{eqnarray*}
Each term is treated separately, integrating by parts
\begin{eqnarray*}
\lefteqn{\int_{\mathbb{R}}\xi \,\overline{u}\,\wedge u_{5}dx =
\int_{\mathbb{R}}\xi \,\wedge(I -
\partial^{2})\overline{u}\,\wedge u_{5}dx = \int_{\mathbb{R}}\xi
\,\wedge \overline{u}\,\wedge u_{5}dx -
\int_{\mathbb{R}}\xi \,\wedge \overline{u}_{2}\,\wedge u_{5}dx }\\
& = & \int_{\mathbb{R}}\partial^{4}\xi \,\wedge \overline{u}\,\wedge
u_{1}dx + \int_{\mathbb{R}}\partial^{3}\xi \,|\wedge u_{1}|^{2}dx -
3\int_{\mathbb{R}}\partial^{2}\xi \,\wedge \overline{u}_{1}\,\wedge
u_{2}dx - 2\int_{\mathbb{R}}\partial\xi \,|\wedge u_{2}|^{2}dx\\
&  & +\int_{\mathbb{R}}\xi \,\wedge \overline{u}_{2}\,\wedge u_{3}dx
- \int_{\mathbb{R}}\partial^{2}\xi \,\wedge \overline{u}_{2}\,\wedge
u_{3}dx - \int_{\mathbb{R}}\partial\xi \,|\wedge u_{3}|^{2}dx +
\int_{\mathbb{R}}\xi \,\wedge \overline{u}_{3}\,\wedge u_{4}dx.
\end{eqnarray*}
The other terms are calculates in a similar way. Then
\begin{eqnarray*}
\lefteqn{2\,i\,Im\int_{\mathbb{R}}\Gamma u\cdot \overline{u}dx}\\
& = & i\,\partial_{t}\int_{\mathbb{R}}\xi\,|u|^{2}dx -
i\int_{\mathbb{R}}\partial_{t}\xi\,|u|^{2}dx - i\,\beta
\int_{\mathbb{R}}\partial^{4}\xi \,\wedge \overline{u}\,\wedge
u_{1}dx - i\,\beta \int_{\mathbb{R}}\partial^{3}\xi \,|\wedge
u_{1}|^{2}dx\\
&  & +\;3\,i\,\beta \int_{\mathbb{R}}\partial^{2}\xi \,\wedge
\overline{u}_{1}\,\wedge u_{2}dx + 2\,i\,\beta
\int_{\mathbb{R}}\partial\xi \,|\wedge u_{2}|^{2}dx - i\,\beta
\int_{\mathbb{R}}\xi \,\wedge
\overline{u}_{2}\,\wedge u_{3}dx\\
&  &  +\;i\,\beta \int_{\mathbb{R}}\partial^{2}\xi \,\wedge
\overline{u}_{2}\, \wedge u_{3}dx + i\,\beta
\int_{\mathbb{R}}\partial\xi \,|\wedge u_{3}|^{2}dx - i\,\beta
\int_{\mathbb{R}}\xi \,\wedge
\overline{u}_{3}\,\wedge u_{4}dx  \\
&  & -\;i\,\beta \int_{\mathbb{R}}\partial^{4}\xi \,\wedge u\,\wedge
\overline{u}_{1}dx - i\,\beta  \int_{\mathbb{R}}\partial^{3}\xi
\,|\wedge u_{1}|^{2}dx + 3\,i\,\beta
\int_{\mathbb{R}}\partial^{2}\xi \,\wedge
u_{1}\,\wedge \overline{u}_{2}dx\\
&  &  +\;2\,i\,\beta \int_{\mathbb{R}}\partial\xi \,|\wedge
u_{2}|^{2}dx - i\,\beta \int_{\mathbb{R}}\xi\,\wedge u_{2}\,\wedge
\overline{u}_{3}dx + i\,\beta\int_{\mathbb{R}}\partial^{2}\xi
\,\wedge u_{2}\,
\wedge\overline{u}_{3}dx\\
&  &  +\;2\,i\,\beta \int_{\mathbb{R}}\partial\xi \,|\wedge
u_{3}|^{2}dx + i\,\beta \int_{\mathbb{R}}\xi\,\wedge
\overline{u}_{3}\,\wedge u_{4}dx + \omega
\int_{\mathbb{R}}\partial^{3}\xi \,\wedge
\overline{u}\,\wedge u_{1}dx\\
&  &  +\;\omega \int_{\mathbb{R}}\partial^{2}\xi\,|\wedge
u_{1}|^{2}dx - 2\,\omega \int_{\mathbb{R}}\partial\xi\,\wedge
\overline{u}_{1}\,\wedge u_{2}dx - \omega \int_{\mathbb{R}}\xi\,|\wedge u_{2}|^{2}dx\\
&  & -\;\omega
\int_{\mathbb{R}}\partial\xi\,\wedge\overline{u}_{2}\,\wedge u_{3}dx
- \omega \int_{\mathbb{R}}\xi \,|\wedge u_{3}|^{2}dx - \omega
\int_{\mathbb{R}}\partial^{3}\xi\,\wedge u\,\wedge
\overline{u}_{1}dx\\
&  &  - \omega \int_{\mathbb{R}}\partial^{2}\xi \,|\wedge
u_{1}|^{2}dx + 2\,\omega \int_{\mathbb{R}}\partial\xi
\,\wedge u_{1}\,\wedge \overline{u}_{2}dx +
\omega \int_{\mathbb{R}}\xi\,|\wedge u_{2}|^{2}dx\\
&  &  +\;\omega \int_{\mathbb{R}}\partial\xi\,\wedge u_{2}\,\wedge
\overline{u}_{3}dx + \omega \int_{\mathbb{R}}\xi \,|\wedge
u_{3}|^{2}dx + i\,\beta \int_{\mathbb{R}}\partial^{2}\xi \,\wedge
\overline{u}\,\wedge u_{1}dx\\
&  &  +\;i\,\beta \int_{\mathbb{R}}\partial\xi \,|\wedge
u_{1}|^{2}dx - i\,\beta \int_{\mathbb{R}}\xi \,\wedge
\overline{u}_{1}\,\wedge u_{2}dx -
i\,\beta\int_{\mathbb{R}}\xi\,\wedge
\overline{u}_{2}\,\wedge u_{3}dx\\
&  & +\;i\,\beta \int_{\mathbb{R}}\partial^{2}\xi \,\wedge u\,\wedge
\overline{u}_{1}dx + i\,\beta \int_{\mathbb{R}}\partial\xi\,|\wedge
u_{1}|^{2}dx - i\,\beta \int_{\mathbb{R}}\xi\,\wedge u_{1}\,\wedge
\overline{u}_{2}dx\\
&  & -\;i\,\beta \int_{\mathbb{R}}\xi \,\wedge u_{2}\,\wedge
\overline{u}_{3}dx - \omega \int_{\mathbb{R}}\partial \xi \,\wedge
\overline{u}\,\wedge u_{1}dx - \omega \int_{\mathbb{R}}\xi\,|\wedge
u_{1}|^{2}dx - \omega
\int_{\mathbb{R}}\xi \,|\wedge u_{2}|^{2}dx\\
&  & +\;\omega \int_{\mathbb{R}}\partial \xi\,\wedge u\,\wedge
\overline{u}_{1}dx + \omega\int_{\mathbb{R}}\xi \,|\wedge
u_{1}|^{2}dx + \omega\int_{\mathbb{R}}\xi\,|\wedge u_{2}|^{2}dx
\end{eqnarray*}
hence
\begin{eqnarray*}
\lefteqn{2\,i\,Im\int_{\mathbb{R}}\Gamma u\cdot \overline{u}dx  =
i\,\partial_{t}\int_{\mathbb{R}}\xi \,|u|^{2}dx -
i\int_{\mathbb{R}}\partial_{t}\xi\,|u|^{2}dx - i\,\beta
\int_{\mathbb{R}}\partial^{4}\xi
\,(|\wedge u|^{2})_{1}dx }\\
&  & -\, 2\,i\,\beta \int_{\mathbb{R}}\partial^{3}\xi \,|\wedge
u_{1}|^{2}dx + 3\,i\,\beta \int_{\mathbb{R}}\partial^{2}\xi
\,(|\wedge u_{1}|^{2})_{1}dx + 4\,i\,\beta
\int_{\mathbb{R}}\partial\xi \,|\wedge
u_{2}|^{2}dx \\
&  &-\,i\,\beta \int_{\mathbb{R}}\xi \,(|\wedge u_{2}|^{2})_{1}dx +
i\,\beta \int_{\mathbb{R}}\partial^{2}\xi \,(|\wedge
u_{2}|^{2})_{1}dx +
3\,i\,\beta \int_{\mathbb{R}}\partial\xi \,|\wedge u_{3}|^{2}dx \\
&  & +\,2\,i\,\omega\,Im\int_{\mathbb{R}}\partial^{3}\xi\,\wedge
\overline{u}\,\wedge u_{1}dx - 4\,i\,\omega \,
Im\int_{\mathbb{R}}\partial\xi \,\wedge \overline{u}_{1}\,\wedge
u_{2}dx \\
&  & -\;2\,i\,\omega \, Im\int_{\mathbb{R}}\partial\xi \,\wedge
\overline{u}_{2}\,\wedge u_{3}dx + i\,\beta
\int_{\mathbb{R}}\partial^{2}\xi \,(|\wedge u|^{2})_{1}dx +
2\,i\,\beta \int_{\mathbb{R}}\partial \xi \,|\wedge u_{1}|^{2}dx
\\
&  & -\;i\,\beta \int_{\mathbb{R}}\xi \,(|\wedge u_{1}|^{2})_{1}dx -
i\,\beta \int_{\mathbb{R}}\xi \,(|\wedge u_{2}|^{2})_{1}dx -
2\,\omega \,Im\int_{\mathbb{R}}\partial \xi\,\wedge
\overline{u}\,\wedge u_{1}dx
\end{eqnarray*}
then, adding similar terms and cutting the letter $i$ we obtain
\begin{eqnarray*}
\lefteqn{2\,Im\int_{\mathbb{R}}\Gamma u\cdot \overline{u}\,dx =
\partial_{t}\int_{\mathbb{R}}\xi \,|u|^{2}dx -
\int_{\mathbb{R}}\partial_{t}\xi\,|u|^{2}dx + \beta
\int_{\mathbb{R}}\partial^{5}\xi \,|\wedge u|^{2}dx -
5\,\beta \int_{\mathbb{R}}\partial^{3}\xi \,|\wedge u_{1}|^{2}dx }\\
&  & +\;6\,\beta \int_{\mathbb{R}}\partial\xi \,|\wedge u_{2}|^{2}dx
- \beta \int_{\mathbb{R}}\partial^{3}\xi\,|\wedge u_{2}|^{2}dx +
3\,\beta \int_{\mathbb{R}}\partial\xi\,|\wedge
u_{3}|^{2}dx \\
&  & +\;2\,\omega\,Im\int_{\mathbb{R}}\partial^{3}\xi\,\wedge
\overline{u}\,\wedge u_{1}dx - 4\,\omega\,
Im\int_{\mathbb{R}}\partial\xi\,\wedge \overline{u}_{1}\,\wedge
u_{2}dx - 2\,\omega\,Im\int_{\mathbb{R}}\partial\xi \,\wedge
\overline{u}_{2}\,\wedge u_{3}dx \\
&  & -\;\beta \int_{\mathbb{R}}\partial^{3}\xi \,|\wedge u|^{2}dx +
3\,\beta \int_{\mathbb{R}}\partial \xi \,|\wedge u_{1}|^{2}dx -
2\,\omega \,Im\int_{\mathbb{R}}\partial \xi \,\wedge
\overline{u}\,\wedge u_{1}dx
\end{eqnarray*}
then
\begin{eqnarray*}
\lefteqn{|\omega |\int_{\mathbb{R}}\partial \xi \,|\wedge
u_{3}|^{2}dx + |\omega |\int_{\mathbb{R}}\partial \xi \,|\wedge
u_{2}|^{2}dx + 2\,|\omega |\int_{\mathbb{R}}\partial \xi \,|\wedge
u_{1}|^{2}dx
 + 2\,|\omega |\int_{\mathbb{R}}\partial \xi \,|\wedge
u_{2}|^{2}\,dx }\\
&  & +\;|\omega |\int_{\mathbb{R}}\partial \xi \,|\wedge u|^{2}dx +
|\omega |\int_{\mathbb{R}}\partial \xi \,|\wedge u_{1}|^{2}dx +
|\omega |\int_{\mathbb{R}}|\partial ^{3}\xi |\,|\wedge
u|^{2}dx\\
& & +\;|\omega |\int_{\mathbb{R}}|\partial ^{3}\xi |\,|\wedge
u_{1}|^{2}dx + \int_{\mathbb{R}}\partial_{t}\xi\,|u|^{2}dx +
2\,Im\int_{\mathbb{R}}\Gamma u\cdot \overline{u}dx\\
& \geq &
\partial_{t}\int_{\mathbb{R}}\xi \,|u|^{2}dx + 3\,\beta
\int_{\mathbb{R}}\partial\xi \,|\wedge u_{3}|^{2}dx -
\beta \int_{\mathbb{R}}\partial^{3}\xi \,|\wedge u_{2}|^{2}dx
+ 6\,\beta \int_{\mathbb{R}}\partial\xi \,|\wedge u_{2}|^{2}dx\\
&  & -\;5\,\beta \int_{\mathbb{R}}\partial^{3}\xi \,|\wedge
u_{1}|^{2}dx + 3\,\beta \int_{\mathbb{R}}\partial \xi \,|\wedge
u_{1}|^{2}dx + \beta \int_{\mathbb{R}}\partial^{5}\xi \,|\wedge
u|^{2}dx - \beta \int_{\mathbb{R}}\partial^{3}\xi \,|\wedge u|^{2}dx
\end{eqnarray*}
where
\begin{eqnarray*}
&  & 3\;|\omega |\int_{\mathbb{R}}\partial \xi \,|\wedge
u_{2}|^{2}dx + |\omega |\int_{\mathbb{R}}[|\partial^{3}\xi | +
3\,\partial \xi ]\,|\wedge u_{1}|^{2}dx \\
&  & +\;|\omega |\int_{\mathbb{R}}[|\partial^{3}\xi | +
\partial \xi + \partial_{t}\xi]\,|\wedge u|^{2}dx +
2\,Im\int_{\mathbb{R}}\Gamma u\cdot \overline{u}dx \\
& \geq &
\partial_{t}\int_{\mathbb{R}}\xi \,|u|^{2}dx +
\int_{\mathbb{R}}[3\,\beta - |\omega |]\,\partial \xi \,|\wedge
u_{3}|^{2}dx - \beta
\int_{\mathbb{R}}\partial^{3}\xi \,|\wedge u_{2}|^{2}dx \\
&  & +\;6\,\beta \int_{\mathbb{R}}\partial\xi \,|\wedge u_{2}|^{2}dx
- 5\,\beta \int_{\mathbb{R}}\partial^{3}\xi \,|\wedge u_{1}|^{2}dx
+\,3\,\beta \int_{\mathbb{R}}\partial \xi \,|\wedge
u_{1}|^{2}dx \\
&  & +\;\beta \int_{\mathbb{R}}\partial^{5}\xi \,|\wedge u|^{2}dx
- \beta \int_{\mathbb{R}}\partial^{3}\xi\,|\wedge u|^{2}dx \\
& \geq &
\partial_{t}\int_{\mathbb{R}}\xi \,|u|^{2}dx + \beta
\int_{\mathbb{R}}[-\partial^{3}\xi +
5\partial \xi ]\,|\wedge u_{2}|^{2}dx \\
&  & +\,\beta \int_{\mathbb{R}}[-5\,\partial^{3}\xi + 3\partial\xi
]\,|\wedge u_{1}|^{2}dx + \beta \int_{\mathbb{R}}[\partial^{3}\xi -
\partial ^{3}\xi ]\,|\wedge u|^{2}dx
\end{eqnarray*}
using \eqref{e203}, $\wedge u_{n}=(I - (I -
\partial^{2}))\wedge u_{n - 2}=\wedge u_{n - 2} - u_{n - 2}$ for $n$
a positive integer and standard estimates we obtain
\begin{eqnarray*}
Im\int_{\mathbb{R}}\Gamma u\cdot \overline{u}\,dx \geq
\partial_{t}\int_{\mathbb{R}}\xi\,|u|^{2}\,dx - c\int_{\mathbb{R}}\xi
\,|u|^{2}\,dx.
\end{eqnarray*}
Multiply this equation by $e^{-Mt},$ and integrate with respect to
$t$ for $t\in [0,\,T]$ and $u\in {\cal D}$
\begin{eqnarray*}
\lefteqn{Im\int_{0}^{T}\int_{\mathbb{R}}e^{-Mt}\,\Gamma u\cdot
\overline{u}\,dx\,dt \geq  \int_{0}^{T}e^{-Mt}\left
(\partial_{t}\int_{\mathbb{R}}\xi\,|u|^{2}dx\right)dt -
c\int_{0}^{T}\int_{\mathbb{R}}\xi\,e^{-Mt}\,|u|^{2}dx\,dt } \\
& = & e^{-Mt}\int_{\mathbb{R}}\xi\,|u|^{2}dx\;\big /_{0}^{T} +
M\int_{0}^{T}\int_{\mathbb{R}}\xi\,e^{-Mt}\,|u|^{2}dx\,dt -
c\int_{0}^{T}\int_{\mathbb{R}}\xi\,e^{-Mt}\,|u|^{2}dx\,dt\\
& = & e^{-Mt}\int_{\mathbb{R}}\xi(x,\,T)\,|u(x,\,T)|^{2}dx +
M\int_{0}^{T}\int_{\mathbb{R}}\xi\,e^{-Mt}\,|u|^{2}dx\,dt -
c\int_{0}^{T}\int_{\mathbb{R}}\xi\,e^{-Mt}\,|u|^{2}dx\,dt.
\end{eqnarray*}
Thus
\begin{eqnarray*}
\lefteqn{<\Gamma
u,\,u>=Im\int_{0}^{T}\int_{\mathbb{R}}e^{-Mt}\,\Gamma u\cdot
\overline{u}\,dx\,dt} \\
& \geq &
 e^{-Mt}\int_{\mathbb{R}}\xi(x,\,T)\,|u(x,\,T)|^{2}dx +
 (M - c)\int_{0}^{T}\int_{\mathbb{R}}\xi\,e^{-Mt}\,|u|^{2}dx\,dt \\
& \geq & \int_{0}^{T}\int_{\mathbb{R}}\xi\,e^{-Mt}\,|u|^{2}dx\,dt
\end{eqnarray*}
provided that $M$ is chosen large enough. Then $<\Gamma u,\,u>\geq
<u,\,u>,$ for all $u\in {\cal D}.$ Let $\Gamma ^{*}$ be the formal
adjoint of $\Gamma $ defined by $\Gamma^{*}=\xi(-i\,\partial_{t} -
i\,\beta\,\wedge \partial^{5} - \omega\,\wedge
\partial^{4} + i\,\beta \,\wedge \partial^{3} + \omega \,\wedge
\partial^{2}).$ Let ${\cal D}^{*}=\{w\in C_{0}^{\infty}(\mathbb{R} \times
[0,\,T]):\;w(x,\,T)=0\,\}.$ In a similar way we prove that
\begin{eqnarray*}
<\Gamma^{*}w,\,w> \;\geq \;<w,\,w>,\quad \forall \;w\in {\cal
D}^{*}.
\end{eqnarray*}
>From this equation, we have that $\Gamma^{*}$ is one-one.
Therefore, $<\Gamma^{*}w,\,\Gamma^{*}v>$ is an inner product on
${\cal D}^{*}.$ We denote by $X$ the completion of ${\cal D}^{*}$
with respect to this inner product. By Riesz's Representation
Theorem, there exists a unique solution $V\in X,$ such that for
any $w\in {\cal D}^{*},$ $<\xi
b^{(1)},\,w>=<\Gamma^{*}V,\,\Gamma^{*}w>$ where we use that $\xi
\,b^{(1)}\in X.$ Then if $v=\Gamma^{*}V$ we have
$<v,\,\Gamma^{*}w>=<\xi b^{(1)},\,w>$ or
$<\Gamma^{*}w,\,v>=<w,\,\xi b^{(1)}>.$ Hence, $v=\Gamma^{*}V$ is a
weak solution of $\Gamma v=\xi b^{(1)}$ with $v\in
L^{2}(\mathbb{R} \times [0,\,T])\simeq
L^{2}([0,\,T]:\,L^{2}(\mathbb{R})).$\\
\\
{\it Remark.} To obtain higher regularity of the solution, we repeat
the proof with higher derivatives. It is a standard approximation
procedure to obtain a result for general initial data.\\
\\
The next step is to estimate the corresponding solutions $v=v(x,\,t)$
of the equation \eqref{e403} via the coefficients of that equation.\\
\\
The following estimate is related to the existence of solutions
theorem.\\
\\
{\bf Lemma 4.2.} {\it Let $|\omega |<3\,\beta $ and $0< \gamma _{1}
\leq \xi \leq \gamma_{2},$ with $\gamma_{2},\,\gamma_{2}$ real
constants. Let $v,\,z\in C^{k}([0,\,+\infty):\;H^{N}(\mathbb{R}))$
for all $k,\,N$ which satisfy \eqref{e403}. For each integer $\alpha
$ there exist positive nondecreasing functions $G$ and $F$ such that
for all $t\geq 0$}
\begin{eqnarray}
\label{e405}\partial_{t}\int_{\mathbb{R}}\xi
\,|v_{\alpha}|^{2}dx\leq G(||z||_{\lambda})\,||v||_{\alpha}^{2} +
F(||z||_{\alpha})
\end{eqnarray}
{\it where $||\;\cdot \;||_{\alpha}$ is the norm in
$H^{\alpha}(\mathbb{R})$
and $\lambda =\max\{1,\,\alpha\}.$}\\
\\
{\it Proof.} Differentiating $\alpha $-times the equation
\eqref{e403}, for some $\alpha \geq 0$ we have
\begin{eqnarray}
\label{e406}-i\;\partial_{t}v_{\alpha} + i\,\beta\wedge v_{\alpha +
5} + \omega \wedge v_{\alpha + 4} - i\,\beta\wedge v_{\alpha + 3} +
\sum_{j=3}^{\alpha + 2}h^{(j)}\wedge v_{j} + (|z|^{2})_{\alpha +
2}\wedge v + p(\wedge z_{\alpha + 1},\,\ldots )=0
\end{eqnarray}
where $h^{(j)}$ is a smooth function depending on $|\wedge
z|^{2},\,\ldots \,$ with $i=2 + \alpha - j.$ For $\alpha\geq 2,$
$p(\wedge z_{\alpha + 1},\,\ldots )$ depends at most linearly on
$\wedge z_{\alpha + 1},$ while for $\alpha =2,$ $p(\wedge z_{\alpha
+ 1},\,\ldots
)$ depends at most quadratically on $\wedge z_{\alpha + 1}.$\\
We multiply equation \eqref{e406} by $\xi\,\overline{v}_{\alpha}$
and integrate over $x\in \mathbb{R}$
\begin{eqnarray*}
\lefteqn{-\,i\int_{\mathbb{R}}\xi
\,\overline{v}_{\alpha}\,\partial_{t}v_{\alpha}dx +
i\,\beta\int_{\mathbb{R}}\xi \,\overline{v}_{\alpha}\wedge v_{\alpha
+ 5}dx + \omega \int_{\mathbb{R}}\xi \,\overline{v}_{\alpha}\,\wedge
v_{\alpha + 4}dx -
i\,\beta\int_{\mathbb{R}}\xi\,\overline{v}_{\alpha}\wedge
v_{\alpha + 3}dx} \\
&  & +\sum_{j=3}^{\alpha + 2}h^{(j)}\int_{\mathbb{R}}\xi
\,\overline{v}_{\alpha}\wedge v_{j}dx +
\int_{\mathbb{R}}\xi\,(|z|^{2})_{\alpha +
2}\overline{v}_{\alpha}\wedge vdx
 +
 \int_{\mathbb{R}}\xi\,\overline{v}_{\alpha}
 p(\wedge z_{\alpha + 1},\,\ldots )dx=0
\end{eqnarray*}
and applying conjugate
\begin{eqnarray*}
\lefteqn{i\int_{\mathbb{R}}\xi
\,v_{\alpha}\,\partial_{t}\overline{v}_{\alpha}dx -
i\,\beta\int_{\mathbb{R}}\xi\,v_{\alpha}\wedge \overline{v}_{\alpha
+ 5}dx + \omega \int_{\mathbb{R}}\xi \,v_{\alpha}\wedge
\overline{v}_{\alpha + 4}dx +
i\,\beta\int_{\mathbb{R}}\xi\,v_{\alpha}\wedge
\overline{v}_{\alpha + 3}dx}\\
&  & +\sum_{j=3}^{\alpha + 2}h^{(j)}\int_{\mathbb{R}}\xi
\,v_{\alpha}\wedge \overline{v}_{j}dx + \int_{\mathbb{R}}\xi
\,(|z|^{2})_{\alpha + 2}v_{\alpha}\wedge \overline{v}dx +
\int_{\mathbb{R}}\xi \,\overline{v}_{\alpha}
 p(\wedge z_{\alpha + 1},\,\ldots )dx=0.
\end{eqnarray*}
Subtracting, it follows that
\begin{eqnarray}
\lefteqn{-\,i\,\partial_{t}\int_{\mathbb{R}}\xi\,|v_{\alpha}|^{2}dx
+ i\,\int_{\mathbb{R}}\partial_{t}\xi\,|v_{\alpha}|^{2}dx +
i\,\beta\int_{\mathbb{R}}\xi\,\overline{v}_{\alpha}\wedge v_{\alpha
+ 5}dx + i\,\beta\int_{\mathbb{R}}\xi
\,v_{\alpha}\wedge\overline{v}_{\alpha + 5}dx } \nonumber \\
&  & +\;\omega \int_{\mathbb{R}}\xi \,\overline{v}_{\alpha}\wedge
v_{\alpha + 4}dx - \omega \int_{\mathbb{R}}\xi \,v_{\alpha}\wedge
\overline{v}_{\alpha + 4}dx - i\,\beta\int_{\mathbb{R}}\xi
\,\overline{v}_{\alpha}\wedge v_{\alpha + 3}dx -
i\,\beta\int_{\mathbb{R}}\xi \,v_{\alpha}\wedge\overline{v}_{\alpha
+ 3}dx \nonumber \\
\label{e407}&  & +\sum_{j=3}^{\alpha + 2}h^{(j)}\int_{\mathbb{R}}\xi
\,\overline{v}_{\alpha}\wedge v_{j}dx - \sum_{j=3}^{\alpha +
2}h^{(j)}\int_{\mathbb{R}}\xi \,v_{\alpha}\wedge\overline{v}_{j}dx +
\int_{\mathbb{R}}\xi
\,(|z|^{2})_{\alpha + 2}v_{\alpha}\wedge\overline{v}dx \\
&  & -\int_{\mathbb{R}}\xi \,(|z|^{2})_{\alpha +
2}\overline{v}_{\alpha}\wedge v\,dx + \int_{\mathbb{R}}\xi
\,\overline{v}_{\alpha}p(\wedge z_{\alpha + 1},\,\ldots )\,dx -
\int_{\mathbb{R}}\xi\,v_{\alpha}\,p(\wedge z_{\alpha + 1},\,\ldots
)dx=0. \nonumber
\end{eqnarray}
Each term is treated separately, integrating by parts
\begin{eqnarray*}
\lefteqn{\int_{\mathbb{R}}\xi\,\overline{v}_{\alpha}\wedge v_{\alpha
+ 5}dx = \int_{\mathbb{R}}\xi\wedge(I -
\partial^{2})\overline{v}_{\alpha}
\wedge v_{\alpha + 5}dx } \\
& = & \int_{\mathbb{R}}\xi\,\wedge\overline{v}_{\alpha}\wedge
v_{\alpha + 5}dx -
\int_{\mathbb{R}}\xi\wedge\overline{v}_{\alpha + 2}\wedge v_{\alpha + 5}dx\\
& = & \int_{\mathbb{R}}\partial^{4}\xi\wedge\overline{v}_{\alpha
}\,\wedge v_{\alpha + 1}dx +
\int_{\mathbb{R}}\partial^{3}\xi\,|\wedge v_{\alpha + 1}|^{2}dx -
3\int_{\mathbb{R}}\partial^{2}\xi\wedge\overline{v}_{\alpha + 1}
\wedge v_{\alpha + 2}dx \\
&  & -\,2\int_{\mathbb{R}}\partial\xi\,|\wedge\overline{v}_{\alpha +
2}|^{2}dx + \int_{\mathbb{R}}\xi\wedge\overline{v}_{\alpha +
2}\,\wedge v_{\alpha + 3}dx -
\int_{\mathbb{R}}\partial^{2}\xi\wedge\overline{v}_{\alpha + 2}
\wedge v_{\alpha + 3}dx \\
&  & -\,2\int_{\mathbb{R}}\partial \xi \,|\wedge v_{\alpha +
3}|^{2}dx - \int_{\mathbb{R}}\xi\,\wedge\overline{v}_{\alpha +
4}\wedge v_{\alpha + 3}dx.
\end{eqnarray*}
The other terms are calculated in a similar way. Hence in
\eqref{e407} we have performing straightforward calculations as
above
\begin{eqnarray*}
\lefteqn{-\,\partial_{t}\int_{\mathbb{R}}\xi\,|v_{\alpha}|^{2}dx +
\int_{\mathbb{R}}\partial_{t}\xi\,|v_{\alpha}|^{2}dx - \beta
\int_{\mathbb{R}}\partial^{5}\xi \,|\wedge v_{\alpha}|^{2}dx +
2\,\beta \int_{\mathbb{R}}\partial^{3}\xi\,
|\wedge v_{\alpha + 1}|^{2}dx } \\
&  & +\,3\,\beta \int_{\mathbb{R}}\partial^{3}\xi\,|\wedge v_{\alpha
+ 1}|^{2}dx - 4\,\beta \int_{\mathbb{R}}\partial\xi\,|\wedge
v_{\alpha + 2}|^{2}dx - \beta
\int_{\mathbb{R}}\partial\xi\,|\wedge v_{\alpha + 2}|^{2}dx\\
&  & +\;\beta \int_{\mathbb{R}}\partial^{2}\xi\,|\wedge v_{\alpha +
2}|^{2}dx - 3\,\beta \int_{\mathbb{R}}\partial\xi\,|\wedge v_{\alpha
+ 3}|^{2}\,dx - 2\,\omega \,Im\int_{\mathbb{R}}\partial^{3}\xi
\,\wedge
\overline{v}_{\alpha}\wedge v_{\alpha + 1}dx \\
&  & +\;4\,\omega \,Im\int_{\mathbb{R}}\partial\xi \wedge
\overline{v}_{\alpha + 1}\wedge v_{\alpha + 2}dx + 2\,\omega
\,Im\int_{\mathbb{R}}\partial\xi\wedge \overline{v}_{\alpha +
2}\wedge v_{\alpha + 3}dx \\
&  & +\;2\,\beta \,Im\int_{\mathbb{R}}\partial \xi
\wedge\overline{v}_{\alpha}\wedge v_{\alpha + 2}dx + 2\,\beta
\,Im\int_{\mathbb{R}}\xi\wedge\overline {v}_{\alpha + 1}
\wedge v_{\alpha + 2}dx \\
&  & -\,\beta \int_{\mathbb{R}}\partial \xi
\,|\wedge\overline{v}_{\alpha + 2}|^{2}dx + 2\,\sum_{j=3}^{\alpha +
2}h^{(j)}\,Im\int_{\mathbb{R}}\xi \,\overline{v}_{\alpha}\wedge
v_{j}dx \\
&  & +\;2\,Im\int_{\mathbb{R}}\xi \,(|z|^{2})_{\alpha +
2}\,v_{\alpha}\wedge\overline{v}dx + 2\,Im\int_{\mathbb{R}}\xi
\,\overline{v}_{\alpha}p(\wedge z_{\alpha + 1},\,\ldots )\,dx=0
\end{eqnarray*}
then
\begin{eqnarray*}
\lefteqn{-\,\partial_{t}\int_{\mathbb{R}}\xi \,|v_{\alpha}|^{2}dx +
\int_{\mathbb{R}}\partial_{t}\xi\,|v_{\alpha}|^{2}dx - 3\,\beta
\int_{\mathbb{R}}\partial\xi\,|\wedge v_{\alpha + 3}|^{2}dx + \beta
\int_{\mathbb{R}}\partial^{2}\xi\,|\wedge v_{\alpha + 2}|^{2}dx } \\
&  & -\,6\,\beta \int_{\mathbb{R}}\partial \xi
\,|\wedge\overline{v}_{\alpha + 2}|^{2}dx + 5\,\beta
\int_{\mathbb{R}}\partial^{3}\xi\,|\wedge v_{\alpha + 1}|^{2}dx -
\beta \int_{\mathbb{R}}\partial^{5}\xi\,|\wedge
v_{\alpha}|^{2}dx \\
& = & -\,2\,\omega\,Im\int_{\mathbb{R}}\partial\xi\,\wedge
\overline{v}_{\alpha + 2}\wedge v_{\alpha + 3}dx - 4\,\omega
\,Im\int_{\mathbb{R}}\partial\xi\wedge \overline{v}_{\alpha +
1}\wedge
v_{\alpha + 2}dx\\
&  & -\,2\,\beta\,Im\int_{\mathbb{R}}\xi\wedge\overline {v}_{\alpha
+ 1}
 \wedge v_{\alpha + 2}dx - 2\,\beta \,Im\int_{\mathbb{R}}\partial \xi
 \wedge\overline{v}_{\alpha}\wedge v_{\alpha + 2}dx\\
&  & +\,2\,\omega\,Im\int_{\mathbb{R}}\partial^{3}\xi \wedge
\overline{v}_{\alpha}\wedge v_{\alpha + 1}dx - 2\,\sum_{j=3}^{\alpha
+ 2}h^{(j)}\,Im\int_{\mathbb{R}}\xi \,
\overline{v}_{\alpha}\wedge v_{j}dx \\
&  & -\,2\,Im\int_{\mathbb{R}}\xi \,(|z|^{2})_{\alpha +
2}v_{\alpha}\wedge\overline{v}dx - 2\,Im\int_{\mathbb{R}}\xi
\,\overline{v}_{\alpha}p(\wedge z_{\alpha + 1},\,\ldots )\,dx
\end{eqnarray*}
hence,
\begin{eqnarray*}
\lefteqn{\partial_{t}\int_{\mathbb{R}}\xi \,|v_{\alpha}|^{2}dx -
\int_{\mathbb{R}}\partial_{t}\xi\,|v_{\alpha}|^{2}dx + 3\,\beta
\int_{\mathbb{R}}\partial\xi\,|\wedge v_{\alpha + 3}|^{2}dx - \beta
\int_{\mathbb{R}}\partial^{2}\xi\,|\wedge v_{\alpha + 2}|^{2}dx } \\
&  & +\, 6\,\beta \int_{\mathbb{R}}\partial\xi
\,|\wedge\overline{v}_{\alpha + 2}|^{2}dx - 5\,\beta
\int_{\mathbb{R}}\partial^{3}\xi\,|\wedge v_{\alpha + 1}|^{2}dx +
\beta
\int_{\mathbb{R}}\partial^{5}\xi \,|\wedge v_{\alpha}|^{2}dx \\
& = & 2\,\omega \,Im\int_{\mathbb{R}}\partial\xi\,\wedge
\overline{v}_{\alpha + 2}\wedge v_{\alpha + 3}dx + 4\,\omega
\,Im\int_{\mathbb{R}}\partial\xi\wedge \overline{v}_{\alpha +
1}\wedge
v_{\alpha + 2}dx\\
&  & +\,2\,\beta \,Im\int_{\mathbb{R}}\xi\wedge\overline {v}_{\alpha
+ 1}
 \wedge v_{\alpha + 2}dx + 2\,\beta \,Im\int_{\mathbb{R}}\partial \xi
 \wedge\overline{v}_{\alpha}\wedge v_{\alpha + 2}dx\\
&  & -\,2\,\omega \,Im\int_{\mathbb{R}}\partial^{3}\xi \wedge
\overline{v}_{\alpha}\wedge v_{\alpha + 1}dx + 2\,\sum_{j=3}^{\alpha
+ 2}h^{(j)}\,Im\int_{\mathbb{R}}\xi \,
\overline{v}_{\alpha}\wedge v_{j}dx \\
&  & +\; 2\,Im\int_{\mathbb{R}}\xi \,(|z|^{2})_{\alpha +
2}v_{\alpha}\wedge\overline{v}dx + 2\,Im\int_{\mathbb{R}}\xi \,
\overline{v}_{\alpha}p(\wedge z_{\alpha + 1},\,\ldots )\,dx\\
& \leq & |\omega |\int_{\mathbb{R}}\partial\xi \,|\wedge v_{\alpha +
2}|^{2}\,dx + |\omega |\int_{\mathbb{R}}\partial\xi |\wedge
v_{\alpha + 3}|^{2}dx + 2\,|\omega |\int_{\mathbb{R}}\partial\xi
\,|\wedge v_{\alpha + 1}|^{2}dx \\
&  & +\;2\,|\omega |\int_{\mathbb{R}}\partial\xi |\wedge v_{\alpha +
2}|^{2}dx + |\beta |\int_{\mathbb{R}}\xi \,|\wedge v_{\alpha +
1}|^{2}dx
+ |\beta |\int_{\mathbb{R}}\xi \,|\wedge v_{\alpha + 2}|^{2}dx\\
&  & +\,|\beta |\int_{\mathbb{R}}\partial \xi |\wedge
v_{\alpha}|^{2}dx + |\beta |\int_{\mathbb{R}}\partial \xi |\wedge
v_{\alpha + 2}|^{2}dx + |\omega
|\int_{\mathbb{R}}\partial^{3}\xi |\wedge v_{\alpha}|^{2}dx\\
&  & +\;|\omega |\int_{\mathbb{R}}\partial^{3}\xi |\wedge v_{\alpha
+ 1}|^{2}dx + 2\,\left |\sum_{j=3}^{\alpha +
2}h^{(j)}\int_{\mathbb{R}}\xi\,\overline{v}_{\alpha}\wedge
v_{j}dx\right | + 2\,\left |\int_{\mathbb{R}}\xi \,(|z|^{2})_{\alpha
+ 2}\,v_{\alpha}\wedge\overline{v}dx\right |
\\
&  & +\;2\,\left |\int_{\mathbb{R}}\xi
\,\overline{v}_{\alpha}p(\wedge z_{\alpha + 1},\,\ldots )\,dx\right
|
\end{eqnarray*}
where
\begin{eqnarray*}
\lefteqn{\partial_{t}\int_{\mathbb{R}}\xi \,|v_{\alpha}|^{2}dx } \\
& \leq &  -\int_{\mathbb{R}}(3\,\beta - |\omega |)\partial \xi
\,|\wedge v_{\alpha + 3}|^{2}dx + \int_{\mathbb{R}}[\beta
\,\partial^{2}\xi - 6\,\beta \,\partial \xi + 3\,|\omega
|\,\partial \xi + |\beta |\,\partial \xi
+ |\beta |\,\xi ]\,|\wedge v_{\alpha + 2}|^{2}dx\\
&  & +\int_{\mathbb{R}}[5\beta\partial^{3}\xi + |\omega
|\partial^{3}\xi + 2\,|\omega |\partial\xi + |\beta |\,\xi
]\,|\wedge v_{\alpha + 1}|^{2}dx  +
\int_{\mathbb{R}}[\partial_{t}\xi + \beta \,\partial^{5}\xi +
|\omega |\,\partial^{3}\xi
+ |\beta |\,\partial \xi ]\,|\wedge v_{\alpha }|^{2}dx\\
&  & +\;2\left |\sum_{j=3}^{\alpha + 2}h^{(j)}\int_{\mathbb{R}}\xi
\,\overline{v}_{\alpha}\wedge v_{j}dx\right | + 2\left
|\int_{\mathbb{R}}\xi \,(|z|^{2})_{\alpha +
2}\,v_{\alpha}\wedge\overline{v}dx\right | + 2\left
|\int_{\mathbb{R}}\xi \,\overline{v}_{\alpha}p(\wedge z_{\alpha +
1},\,\ldots )\,dx\right |.
\end{eqnarray*}
using that $|\omega |<3\,\beta $ we have that the first term in
the right hand side of the above expression is not positive.
Hence,
\begin{eqnarray*}
\lefteqn{\partial_{t}\int_{\mathbb{R}}\xi \,|v_{\alpha}|^{2}dx } \\
& \leq & \int_{\mathbb{R}}[\beta \,\partial^{2}\xi - 6\,\beta
\,\partial \xi + 3\,|\omega |\,\partial \xi + |\beta |\,\partial
\xi + |\beta |\,\xi ]\,|\wedge v_{\alpha + 2}|^{2}dx  \\
&  & +\int_{\mathbb{R}}[5\,\beta \,\partial^{3}\xi + |\omega
|\,\partial^{3}\xi + 2\,|\omega |\,\partial \xi + |\beta |\,\xi
]\,|\wedge v_{\alpha + 1}|^{2}dx +
\int_{\mathbb{R}}[\partial_{t}\xi + \beta \,\partial^{5}\xi +
|\omega |\,\partial^{3}\xi
+ |\beta |\,\partial \xi ]\,|\wedge v_{\alpha }|^{2}dx\\
&  & +\,2\left |\sum_{j=3}^{\alpha + 2}h^{(j)}\int_{\mathbb{R}}\xi
\,\overline{v}_{\alpha}\wedge v_{j}dx\right | + 2\left
|\int_{\mathbb{R}}\xi \,(|z|^{2})_{\alpha +
2}\,v_{\alpha}\wedge\overline{v}\,dx\right | + 2\left
|\int_{\mathbb{R}}\xi\,\overline{v}_{\alpha}p(\wedge z_{\alpha +
1},\,\ldots)\,dx\right |.
\end{eqnarray*}
Using that $\wedge v_{n} =\wedge v_{n - 2} - v_{n - 2}$ and a
standard estimate, the lemma follows.
\renewcommand{\theequation}{\thesection.\arabic{equation}}
\setcounter{equation}{0}\section{Uniqueness and Existence of a Local
Solution} In this section, we study the uniqueness and the existence
of local strong solutions in the Sobolev space $H^{N}(\mathbb{R})$
for $N\geq 3$ for the problem \eqref{e204}. To establish the
existence of strong solutions for
\eqref{e204} we use the a priori estimate together with an
approximation procedure.\\
\\
{\bf Theorem 5.1}(Uniqueness). {\it Let $|\omega |<3\,\beta ,$
$u_{0}(x)\in H^{N}(\mathbb{R})$ with $N\geq 3$ and $0<T<+\infty .$
Then there is at most one strong solution $u\in
L^{\infty}([0,\,T]:\,H^{N}(\mathbb{R}))$ of \eqref{e204} with
initial data
$u(x,\,0)=u_{0}(x).$}\\
\\
{\it Proof.} Assume that $u,\,v\in
L^{\infty}([0,\,T]:\,H^{N}(\mathbb{R}))$ are two solutions of
\eqref{e204} with $u_{t},$ $v_{t}$ $\in L^{\infty}([0,\,T]:\,H^{N -
3}(\mathbb{R})),$ and with the same initial data. Then
\begin{eqnarray}
\label{e501}i\,(u - v)_{t} + i\,\beta \,(u - v)_{3} + \omega \,(u -
v)_{2} + |u|^{2}\,u - |v|^{2}\,v = 0
\end{eqnarray}
with $(u - v)(x,\,0)=0.$ By \eqref{e501}
\begin{eqnarray*}
i\,(u - v)_{t} + i\,\beta \,(u - v)_{3} +
\omega \,(u - v)_{2} + |u|^{2}\,(u - v) + (|u|^{2} - |v|^{2})\,v = 0
\end{eqnarray*}
or
\begin{eqnarray}
\label{e502}i\,(u - v)_{t} + i\,\beta \,(u - v)_{3} + \omega \,(u -
v)_{2} + |u|^{2}\,(u - v) + (|u| - |v|)\,(|u| + |v|)\,v = 0.
\end{eqnarray}
Multiplying \eqref{e502} by $\xi\overline{(u - v)}$ we have
\begin{eqnarray*}
\lefteqn{i\,\xi \,\overline{(u - v)}\,(u - v)_{t} + i\,\beta \,\xi
\, \overline{(u - v)}\,(u - v)_{3} + \alpha \,\xi \,\overline{(u -
v)}\,
(u - v)_{2} } \\
&  & +\,|u|^{2}\,|u - v|^{2} +
\xi \,\overline{(u - v)}\,(|u| - |v|)\,(|u| + |v|)\,v = 0.\\
\\
\lefteqn{-\,i\,\xi \,(u - v)\,\overline{(u - v)}_{t} - i\,\beta \,
\xi \,(u - v)\,\overline{(u - v)}_{3} + \alpha \,\xi \,
(u - v)\,\overline{(u - v)}_{2} } \\
&  & +\,|u|^{2}\,|u - v|^{2} + \xi \,(u - v)\,(|u| - |v|)\,(|u| +
|v|)\,\overline{v} = 0.\quad (\mbox{applying conjugate)}
\end{eqnarray*}
Subtracting and integrating over $x\in \mathbb{R}$ we obtain
\begin{eqnarray}
\lefteqn{i\,\partial_{t}\int_{\mathbb{R}}\xi \,|u - v|^{2}dx -
i\int_{\mathbb{R}}\partial_{t}\xi\,|u - v|^{2}dx + i\,\beta
\int_{\mathbb{R}}\xi \,\overline{(u - v)}\,(u - v)_{3}dx } \nonumber \\
&  & +\;i\,\beta \int_{\mathbb{R}}\xi \,(u - v)\,\overline{(u -
v)}_{3}dx +\,\omega \int_{\mathbb{R}}\xi \,\overline{(u - v)}\,(u -
v)_{2}dx \nonumber \\
\label{e503}&  & -\,\omega \int_{\mathbb{R}}\xi \,(u -
v)\,\overline{(u - v)}_{2}dx + 2\,i\,Im\int_{\mathbb{R}}\xi
\,\overline{(u - v)}\,(|u| - |v|)\,(|u| + |v|)\,v\,dx = 0\quad
\end{eqnarray}
Each term is treated separately, integrating by parts
\begin{eqnarray*}
\lefteqn{\int_{\mathbb{R}}\xi \,\overline{(u - v)}\,(u - v)_{3}dx} \\
& = & \int_{\mathbb{R}}\partial^{2}\xi\,\overline{(u - v)}\,(u -
v)_{1}dx + 2\int_{\mathbb{R}}\partial \xi \,|(u - v)_{1}|^{2}dx +
\int_{\mathbb{R}}\xi \,(u - v)_{1}\,\overline{(u - v)}_{2}dx.
\end{eqnarray*}
The other terms are calculated in a similar way. Hence in
\eqref{e503} we have
\begin{eqnarray*}
\lefteqn{i\,\partial_{t}\int_{\mathbb{R}}\xi \,|u - v|^{2}dx -
i\int_{\mathbb{R}}\partial_{t}\xi\,|u - v|^{2}dx + i\,\beta
\int_{\mathbb{R}}\partial^{2}\xi \,\overline{(u - v)}\,(u - v)_{1}dx } \\
&  & +\,2\,i\,\beta \int_{\mathbb{R}}\partial \xi \,|(u -
v)_{1}|^{2}dx + i\,\beta \int_{\mathbb{R}}\xi\,(u -
v)_{1}\,\overline{(u - v)}_{2}dx + i\,\beta
\int_{\mathbb{R}}\partial^{2}\xi\,(u - v)\,\overline{(u -
v)}_{1}dx \\
&  & +\,i\,\beta \int_{\mathbb{R}}\partial \xi\,|(u - v)_{1}|^{2}dx
- i\,\beta \int_{\mathbb{R}}\xi\,(u - v)_{1}\,\overline{(u -
v)}_{2}dx - \omega \int_{\mathbb{R}}\partial \xi \,\overline{(u -
v)}\,(u - v)_{1}dx
\\
&  & -\,\omega \int_{\mathbb{R}}\xi\,|(u - v)_{1}|^{2}dx + \omega
\int_{\mathbb{R}}\partial \xi \,(u - v)\,\overline{(u - v)}_{1}dx
+ \omega
\int_{\mathbb{R}}\xi\,|(u - v)_{1}|^{2}dx \\
&  & +\,2\,i\,Im\int_{\mathbb{R}}\xi \,\overline{(u - v)}\,(|u| -
|v|)\,(|u| + |v|)\,v\,dx = 0
\end{eqnarray*}
then
\begin{eqnarray*}
\lefteqn{i\,\partial_{t}\int_{\mathbb{R}}\xi \,|u - v|^{2}dx -
i\int_{\mathbb{R}}\partial_{t}\xi\,|u - v|^{2}dx + i\,\beta
\int_{\mathbb{R}}\partial^{2}\xi \,(|u - v|^{2})_{1}dx + 3\,i\,\beta
\int_{\mathbb{R}}\partial \xi \,|(u -
v)_{1}|^{2}dx } \\
&  & -\;2\,i\,\omega \,Im\int_{\mathbb{R}}\partial \xi \,
\overline{(u - v)}\,(u - v)_{1}dx + 2\,i\,Im\int_{\mathbb{R}}\xi \,
\overline{(u - v)}\,(|u| - |v|)\,(|u| + |v|)\,v\,dx = 0
\end{eqnarray*}
if and only if
\begin{eqnarray*}
\lefteqn{\partial_{t}\int_{\mathbb{R}}\xi \,|u - v|^{2}dx -
\int_{\mathbb{R}}\partial_{t}\xi\,|u - v|^{2}\,dx + \beta
\int_{\mathbb{R}}\partial^{2}\xi\,(|u - v|^{2})_{1}dx + 3\,\beta
\int_{\mathbb{R}}
\partial\xi\,|(u - v)_{1}|^{2}dx } \\
& = & 2\,\omega\,Im\int_{\mathbb{R}}\partial\xi\,\overline{(u -
v)}\, (u - v)_{1}dx - 2\,Im\int_{\mathbb{R}}\xi\,\overline{(u -
v)}\, (|u| - |v|)\,(|u| + |v|)\,v\,dx \\
& \leq & |\omega|\int_{\mathbb{R}}\partial\xi\,|u - v|^{2}dx +
|\omega |\int_{\mathbb{R}}\partial\xi\,|(u - v)_{1}|^{2}dx +
2\int_{\mathbb{R}}\xi \,|u - v|\,|\;|u| - |v|\;|\,(|u| +
|v|)\,|v|\,dx.
\end{eqnarray*}
Using that $|\;|u| - |v|\;|\leq |u - v|,$ \eqref{e203}
 and standard estimates, we have
\begin{eqnarray*}
&  & \partial_{t}\int_{\mathbb{R}}\xi \,|u - v|^{2}dx +
\int_{\mathbb{R}}[3\,\beta - |\omega|\,]\,\partial \xi \,|(u -
v)_{1}|^{2}dx \leq c\int_{\mathbb{R}}\xi \,|u - v|^{2}dx.
\end{eqnarray*}
Integrating in $t\in [0,\,T],$ using the fact that $(u - v)$
vanishes at $t=0$ and Gronwall's inequality it follows that $u=v.$
This proves the
uniqueness of the solution.\\
\\
We construct the mapping ${\cal
Z}:L^{\infty}([0,\,T]:\,H^{s}(\mathbb{R})) \longmapsto
L^{\infty}([0,\,T]:\,H^{s}(\mathbb{R}))$ where the initial
condition is given by $u^{(n)}(x,\,0)=u_{0}(x)$ and the first
approximation is given by
\begin{eqnarray*}
u^{(0)} & = & u_{0}(x)\\
u^{(n)} & = & {\cal Z}(u^{(n - 1)})\qquad n\geq 1,
\end{eqnarray*}
where $u^{(n - 1)}$ is in place of $z$ in equation \eqref{e403} and
$u^{(n)}$ is in place of $v$ which is the solution of equation
\eqref{e403}. That is
\begin{eqnarray*}
\lefteqn{-\,i\,u_{t}^{(n)} + i\,\beta \,\wedge u_{5}^{(n)} + \omega
\,\wedge u_{4}^{(n)} + (|\wedge u^{(n - 1)}|^{2})_{2}\wedge u^{(n)}
+ 2\,(|\wedge u^{(n - 1)}|^{2})_{1}\wedge
u_{1}^{(n)} } \\
&  & +\;|\wedge u^{(n - 1)}|^{2}\,\wedge u_{2}^{(n)} - (i\,\beta
\,\wedge u_{3}^{(n)} + \omega \,\wedge u_{2}^{(n)} + |\wedge u^{(n -
1)}|^{2}\,\wedge u^{(n)})=0.
\end{eqnarray*}
By Lemma 4.1, $u^{(n)}$ exists and is unique in $C((0,\,+\infty
):\,H^{N}(\mathbb{R})).$ A choice of $c_{0}$ and the use of the a
priori estimate in Section 4 shows that ${\cal
Z}:\mathbb{B}_{c_{0}}(0)\longmapsto \mathbb{B}_{c_{0}}(0)$ where
$\mathbb{B}_{c_{0}}(0)$ is a bounded ball in
$L^{\infty}([0,\,T]:\,H^{s}(\mathbb{R})).$\\
\\
{\bf Theorem 5.2}(Local solution). {\it Let $|\omega |<3\,\beta$
and $N$ an integer $\geq 3.$ If $u_{0}(x) \in H^{N}(\mathbb{R}),$
then there is $T>0$ and $u$ such that $u$ is a strong solution of
\eqref{e204},
$u\in L^{\infty}([0,\,T]:\,H^{N}(\mathbb{R}))$ and $u(x,\,0)=u_{0}(x).$}\\
\\
{\it Proof.} We prove that for $u_{0}(x)\in
H^{\infty}(\mathbb{R})=\bigcap_{k\geq 0}H^{k}(\mathbb{R})$ there
exists a solution $u\in L^{\infty}([0,\,T]:\,H^{N}(\mathbb{R}))$
with initial data $u(x,\,0)=u_{0}(x)$ where the time of existence
$T>0$ only depends on the norm of $u_{0}(x).$ We define a sequence
of approximations to equation \eqref{e403} as
\begin{eqnarray}
i\,v_{t}^{(n)} & = & i\,\beta \wedge v_{5}^{(n)} + \omega \,\wedge
v_{4}^{(n)} - i\,\beta \wedge v_{3}^{(n)}
- \omega \wedge v_{2}^{(n)} + |\wedge v^{(n - 1)}|^{2}\,\wedge v_{2}^{(n)} \nonumber \\
\label{e504}&  & +\;O[\,(|\wedge v^{(n - 1)}|^{2})_{2},\,(|\wedge
v^{(n - 1)}|^{2})_{1},\,\ldots)\,]
\end{eqnarray}
where the initial condition is $v^{(n)}(x,\,0) = u_{0}(x) -
\partial^{2}u_{0}(x).$ The first approximation is given by
$v^{(0)}(x,\,0) = u_{0}(x) - \partial^{2}u_{0}(x).$ Equation
\eqref{e504} is a linear equation at each iteration which can be
solved in any interval of time in which the coefficients are
defined. This is shown in Lemma 4.1. By Lemma 4.2, it follows that
\begin{eqnarray}
\label{e505}\partial_{t}\int_{\mathbb{R}}\xi
\,|v_{\alpha}^{(n)}|^{2}dx\leq G(||v^{(n -
1)}||_{\lambda})\,||v^{(n)}||_{\alpha}^{2} + F(||v^{(n -
1)}||_{\alpha}).
\end{eqnarray}
Choose $\alpha =1$ and let $c\geq ||u_{0} -
\partial^{2}u_{0}||_{1}\geq ||u_{0}||_{3}.$
For each iterate $n,$ $||v^{(n)}(\,\cdot\,,\,t)||$ is continuous
in $t\in [0,\,T]$ and $||v^{(n)}(\,\cdot\,,\,0)||<c.$ Define
$c_{0}=\frac{\gamma_{2}}{2\,\gamma_{1}}\,c^{2} + 1.$ Let
$T_{0}^{(n)}$ be the maximum time such that $||v^{(k)}(\,\cdot
\,,\,t)||_{1}\leq c_{3}$ for $0\leq t\leq T_{0}^{(n)},$ $0\leq
k\leq n.$ Integrating \eqref{e505} over $[0,\,t]$ we have that for
$0\leq t\leq T_{0}^{(n)}$ and $j=0,\,1$
\begin{eqnarray*}
\int_{0}^{t}\left(\partial_{s}\int_{\mathbb{R}}\xi
\,|v_{j}^{(n)}|^{2}dx\right)ds\leq \int_{0}^{t}G\left(||v^{(n -
1)}||_{1}\right)||v^{(n)}||_{j}^{2}ds + \int_{0}^{t}F\left(||v^{(n -
1)}||_{j}\right)ds.
\end{eqnarray*}
It follows that
\begin{eqnarray*}
\int_{\mathbb{R}}\xi(x,\,t)|v_{j}^{(n)}(x,\,t)|^{2}dx & \leq &
\int_{\mathbb{R}}\xi(x,\,0)|v_{j}^{(n)}(x,\,0)|^{2}dx +
\int_{0}^{t}G\left(||v^{(n - 1)}||_{1}\right)||v^{(n)}||_{j}^{2}ds
\\
&  & + \int_{0}^{t}F\left(||v^{(n - 1)}||_{j}\right)ds
\end{eqnarray*}
hence
\begin{eqnarray*}
\gamma_{1}\int_{\mathbb{R}}|v_{j}^{(n)}(x,\,t)|^{2}dx & \leq &
\int_{\mathbb{R}}\xi(x,\,t)|v_{j}^{(n)}(x,\,t)|^{2}dx \\
& \leq & \int_{\mathbb{R}}\xi(x,\,0)|v_{j}^{(n)}(x,\,0)|^{2}dx +
\int_{0}^{t}G\left(||v^{(n - 1)}||_{1}\right)||v^{(n)}||_{j}^{2}ds
\\
&  & + \int_{0}^{t}F\left(||v^{(n - 1)}||_{j}\right)ds
\end{eqnarray*}
and
\begin{eqnarray*}
\int_{\mathbb{R}}|v_{j}^{(n)}|^{2}dx \leq
\frac{\gamma_{2}}{\gamma_{1}}\int_{\mathbb{R}}|v_{j}^{(n)}(x,\,0)|^{2}dx
+ \frac{G(c_{3})}{\gamma_{1}}\,c_{3}^{2}\,t +
\frac{F(c_{3})}{\gamma_{1}}\,t
\end{eqnarray*}
and we obtain for $j=0,\,1$ that
\begin{eqnarray*}
||v^{(n)}||_{1}\leq \frac{\gamma_{2}}{\gamma_{1}}\,c^{2} +
\frac{G(c_{0})}{\gamma_{1}}\,c_{0}^{2}\,t + \frac{F(c_{0})}{\gamma_{1}}\,t.
\end{eqnarray*}
{\it Claim.} $T_{0}^{(n)}$ does not approach to $0.$\\
On the contrary, assume that $T_{0}^{(n)}\rightarrow 0.$ Since
$||v^{(n)}(\,\cdot \,,\,t)||$ is continuous for $t\geq 0,$ there
exists $\tau \in [0,\,T]$ such that
$||v^{(k)}(\,\cdot\,,\,t)||_{1}=c_{0}$ for $0\leq \tau \leq
T_{0}^{(n)},$ $0\leq k\leq n.$ Then
\begin{eqnarray*}
c_{0}^{2}\leq \frac{\gamma_{2}}{\gamma_{1}}\,c^{2} +
\frac{G(c_{0})}{\gamma_{1}}\,c_{0}^{2}\,T_{0}^{(n)} +
\frac{F(c_{0})}{\gamma_{1}}\,T_{0}^{(n)}
\end{eqnarray*}
as $n\rightarrow \infty ,$ we have
\begin{eqnarray*}
\left(\frac{\gamma_{2}}{2\,\gamma_{1}}\,c^{2} + 1\right)^{2} \leq
\frac{\gamma_{2}}{\gamma_{1}}\,c^{2}\quad\mbox{then}\quad
\frac{\gamma_{2}^{2}}{4\,\gamma_{1}^{2}}\,c^{4} + 1\leq 0
\end{eqnarray*}
which is a contradiction. Consequently $T_{0}^{(n)}\not \rightarrow 0.$
Choosing $T=T(c)$ sufficiently small, and $T$ not depending on $n,$ one
concludes that
\begin{eqnarray}
\label{e506}||v^{(n)}||_{1}\leq C
\end{eqnarray}
for $0\leq t\leq T.$ This shows that $T_{0}^{(n)}\geq T.$ Hence,
from \eqref{e506} we imply that there exists a subsequence
$v^{(n_{j})}\equiv v^{(n)}$ such that
\begin{eqnarray}
\label{e507}v^{(n)}\stackrel{*}\rightharpoonup v\quad \mbox{weakly
on}\quad L^{\infty}([0,\,T]:\,H^{1}(\mathbb{R})).
\end{eqnarray}
{\it Claim.} $u=\wedge v$ is a solution.\\
\\
In the linearized equation \eqref{e504} we have
\begin{eqnarray*}
\wedge v_{5}^{(n)}=\wedge(I - (I - \partial^{2}))v_{3}^{(n)} =
\wedge v_{3}^{(n)} - v_{3}^{(n)} =
\partial^{2}(\underbrace{\wedge v_{1}^{(n)}}_{\in L^{2}(\mathbb{R})}) -
\underbrace{ \partial^{2}(v_{1}^{(n)})}_{\in H^{-2}(\mathbb{R})}\in
H^{-2}(\mathbb{R}).
\end{eqnarray*}
Since $\;\wedge= (I - \partial^{2})^{-1}\;$ is bounded in
$\;H^{1}(\mathbb{R}),\;$  $\;\wedge v_{5}^{(n)}\;$ belongs to
$\;H^{-2}(\mathbb{R}).\;$ $\;v^{(n)}\;$ is still bounded in
$L^{\infty}([0,\,T]:\;H^{1}(\mathbb{R}))\hookrightarrow
L^{2}([0,\,T]:\;H^{1}(\mathbb{R}))$ and since
$\wedge:L^{2}(\mathbb{R})\rightarrow H^{2}(\mathbb{R})$ is a bounded
operator,
\begin{eqnarray*}
||\wedge v_{1}^{(n)}||_{H^{2}(\mathbb{R})}\leq
c\,||v_{1}^{(n)}||_{L^{2}(\mathbb{R})}\leq c\,||
v_{1}^{(n)}||_{H^{1}(\mathbb{R})}.
\end{eqnarray*}
Consequently, $\wedge v_{1}^{(n)}$ is bounded in
$L^{2}([0,\,T]:\,H^{2}(\mathbb{R}))\hookrightarrow
L^{2}([0,\,T]:\,L^{2}(\mathbb{R})).$ It follows that
$\partial^{2}(\wedge v_{1}^{(n)})$ is bounded in
$L^{2}([0,\,T]:\,H^{-2}(\mathbb{R})),$ and
\begin{eqnarray}
\label{e508}\wedge v_{5}^{(n)}\quad \mbox{is bounded in}\quad
L^{2}([0,\,T]:\;H^{-2}(\mathbb{R})).
\end{eqnarray}
Similarly, the other terms are bounded. By \eqref{e504},
$v_{t}^{(n)}$ is a sum of terms each of which is the product of a
coefficient, uniformly bounded on $n$ and a function in
$L^{2}([0,\,T]:\,H^{-2}(\mathbb{R}))$ uniformly bounded on $n$
such that $v_{t}^{(n)}$ is bounded in
$L^{2}([0,\,T]:\,H^{-2}(\mathbb{R})).$ On the other hand,
$H_{loc}^{1}(\mathbb{R})\stackrel{c}\hookrightarrow
H_{loc}^{1/2}(\mathbb{R})\hookrightarrow H^{-4}(\mathbb{R}).$ By
Lions-Aubin's compactness Theorem \cite{li1} there is a
subsequence $v^{(n_{j})}\equiv v^{(n)}$ such that
$v^{(n)}\rightarrow v$ strongly on
$L^{2}([0,\,T]:\,H_{loc}^{1/2}(\mathbb{R})).$ Hence, for a
subsequence $v^{(n_{j})}\equiv v^{(n)},$ we have
$v^{(n)}\rightarrow v$ a. e. in
$L^{2}([0,\,T]:\,H_{loc}^{1/2}(\mathbb{R})).$ Moreover, from
\eqref{e508}, $\wedge v_{5}^{(n)}\rightharpoonup \wedge v_{5}$
weakly in $L^{2}([0,\,T]:\,H^{-2}(\mathbb{R})).$ Similarly,
$\wedge v_{2}^{(n)}\rightharpoonup \wedge v_{2}$ weakly in
$L^{2}([0,\,T]:\,H^{-2}(\mathbb{R})).$ Since $||\wedge
v^{(n)}||_{H^{2}(\mathbb{R})}\leq c\,||
v^{(n)}||_{L^{2}(\mathbb{R})}\leq c\,||
v^{(n)}||_{H^{1}(\mathbb{R})}\leq c\,||
v^{(n)}||_{H^{1/2}(\mathbb{R})}$ and $v^{(n)}\rightarrow v$
strongly on $L^{2}([0,\,T]:\,H_{loc}^{1/2}(\mathbb{R}))$ then
$\wedge v^{(n)}\rightarrow \wedge v$ strongly in
$L^{2}([0,\,T]:\,H_{loc}^{2}(\mathbb{R})).$ Thus, the fifth term
on the right hand side of \eqref{e504}, $|\wedge v^{(n -
1)}|^{2}\,\wedge v_{2}^{(n)}\rightharpoonup |\wedge v|^{2}\,\wedge
v_{2}$ weakly in $L^{2}([0,\,T]:\,L_{loc}^{1}(\mathbb{R}))$ as
$\wedge v_{2}^{(n)}\rightharpoonup \wedge v_{2}$ weakly in
$L^{2}([0,\,T]:\,H^{-2}(\mathbb{R}))$ and $|\wedge v^{(n -
1)}|^{2}\rightarrow |\wedge v|^{2}$ strongly on
$L^{2}([0,\,T]:\,H_{loc}^{2}(\mathbb{R})).$ Similarly, the other
terms in \eqref{e504} converge to their limits, implying
$v_{t}^{(n)}\rightharpoonup v_{t}$ weakly in
$L^{2}([0,\,T]:\,L_{loc}^{1}(\mathbb{R})).$ Passing to the limit
\begin{eqnarray*}
i\,v_{t} & = & \partial^{2}(i\,\beta \,\wedge v_{3} + \omega \,\wedge v_{2} +
|\wedge v|^{2}\,\wedge v) - (i\,\beta \,\wedge v_{3} +
\omega \,\wedge v_{2} + |\wedge v|^{2}\,\wedge v)\\
& = & -(I - \partial^{2})(i\,\beta \,\wedge v_{3} + \omega
\,\wedge v_{2} + |\wedge v|^{2}\,\wedge v).
\end{eqnarray*}
Thus $i\,v_{t} + (I - \partial^{2})(i\,\beta \,\wedge v_{3} +
\omega \,\wedge v_{2} + |\wedge v|^{2}\,\wedge v)=0.$ This way, we
have \eqref{e204} for $u=\wedge v.$\\
\\
Now, we prove that there exists a solution of \eqref{e204} with
$u\in L^{\infty}([0,\,T]:\,H^{N}(\mathbb{R}))$ and $N\geq 4,$
where $T$ depends only on the norm of $u_{0}$ in
$H^{3}(\mathbb{R}).$ We already know that there is a solution
$u\in L^{\infty}([0,\,T]:\,H^{3}(\mathbb{R})).$ It is suffices to
show that the approximating sequence $v^{(n)}$ is bounded in
$L^{\infty}([0,\,T]:\,H^{N - 2}(\mathbb{R})).$ Taking $\alpha = N
- 2$ and considering \eqref{e505} for $\alpha \geq 2,$ we define
$c_{N - 2}=\frac{\gamma_{2}}{2\,\gamma_{1}}\,||u_{0}(\cdot)||_{N}
+ 1.$ Let $T_{N - 3}^{(n)}$ be the largest time such that
$||v^{(k)}(\,\cdot\,,\,t)||_{\alpha}\leq c_{N - 3}$ for $0\leq
t\leq T_{N - 3}^{(n)},$ $0\leq k\leq n.$ Integrating \eqref{e505}
over $[0,\,t],$ for $0\leq t\leq T_{N - 3}^{(n)},$ we have
\begin{eqnarray*}
\int_{0}^{t}\left(\partial_{s}\int_{\mathbb{R}}\xi
\,|v_{\alpha}^{(n)}|^{2}dx\right)ds \leq \int_{0}^{t}G\left(||v^{(n
- 1)}||_{\alpha}\right)||v^{(n)}||_{\alpha}^{2}ds +
\int_{0}^{t}F\left(||v^{(n - 1)}||_{\alpha}\right)ds.
\end{eqnarray*}
It follows that
\begin{eqnarray*}
\int_{\mathbb{R}}\xi (x,\,t)\,|v_{\alpha}^{(n)}|^{2}dx & \leq &
\int_{\mathbb{R}}\xi (x,\,0)\, |v_{\alpha}^{(n)}(x,\,0)|^{2}dx +
\int_{0}^{t}G\left(||v^{(n - 1)}||_{\alpha}\right)
||v^{(n)}||_{\alpha}^{2}ds \\
&  & + \int_{0}^{t}F\left(||v^{(n - 1)}||_{\alpha}\right)ds
\end{eqnarray*}
hence
\begin{eqnarray*}
\gamma_{1}\int_{\mathbb{R}}|v_{\alpha}^{(n)}|^{2}dx\leq
\int_{\mathbb{R}}\xi \,|v_{\alpha}^{(n)}|^{2}dx & \leq &
\int_{\mathbb{R}}\xi (x,\,0)\, |v_{\alpha}^{(n)}(x,\,0)|^{2}dx +
\int_{0}^{t}G
\left(||v^{(n - 1)}||_{\alpha}\right)||v^{(n)}||_{\alpha}^{2}ds \\
&  & +\int_{0}^{t}F\left(||v^{(n - 1)}||_{\alpha}\right)ds
\end{eqnarray*}
then
\begin{eqnarray*}
\lefteqn{\int_{\mathbb{R}}|v_{\alpha}^{(n)}|^{2}dx\leq
\frac{\gamma_{2}}{\gamma_{1}}
\int_{\mathbb{R}}|v_{\alpha}^{(n)}(x,\,0)|^{2}dx + \frac{G(c_{N -
3})}{\gamma_{1}}\,c_{N - 3}^{2}\,t +
\frac{F(c_{N - 3})}{\gamma_{1}}\,t } \\
& \leq & \frac{\gamma_{2}}{\gamma_{1}}\,||v_{\alpha}^{(n)}(x,\,0)||_{\alpha}^{2}
+ \frac{G(c_{N - 3})}{\gamma_{1}}\,c_{N - 3}^{2}\,t +
\frac{F(c_{N - 3})}{\gamma_{1}}\,t\\
& \leq & \frac{\gamma_{2}}{\gamma_{1}}\,||u(x,\,0)||_{N}^{2} +
\frac{G(c_{N - 3})}{\gamma_{1}}\,c_{N - 3}^{2}\,t +
\frac{F(c_{N - 3})}{\gamma_{1}}\,t
\end{eqnarray*}
and we obtain
\begin{eqnarray*}
||v_{\alpha}^{(n)}(\,\cdot \,,\,t)||_{\alpha}^{2}dx\leq
\frac{\gamma_{2}}{\gamma_{1}}\,||u(x,\,0)||_{N}^{2} + \frac{G(c_{N -
3})}{\gamma_{1}}\,c_{N - 3}^{2}\,t + \frac{F(c_{N -
3})}{\gamma_{1}}\,t
\end{eqnarray*}
{\it Claim.} $T_{N - 3}^{(n)}$ does not approach to $0.$\\
On the contrary, assume that $T_{N - 3}^{(n)}\rightarrow 0.$ Since
$||v^{(n)}(\,\cdot\,,\,t)||$ is continuous for $t\geq 0,$ there
exists $\tau \in [0,\,T_{N - 3}]$ such that $||v^{(k)}(\,\cdot
\,,\,\tau)||_{\alpha}=c_{N - 3}$ for $0\leq \tau \leq T^{(n)},$
$0\leq k\leq n.$ Then
\begin{eqnarray*}
c_{N - 3}^{2}\leq \frac{\gamma_{2}}{\gamma_{1}}\,||u(x,\,0)||_{N}^{2} +
\frac{G(c_{N - 3})}{\gamma_{1}}\,c_{N - 3}^{2}\,T_{N - 3}^{(n)} +
\frac{F(c_{N - 3})}{\gamma_{1}}\,T_{N - 3}^{(n)}
\end{eqnarray*}
as $n\rightarrow +\infty ,$ and we have
\begin{eqnarray*}
\left(\frac{\gamma_{2}}{2\,\gamma_{1}}\,||u(x,\,0)||_{N}^{2} + 1\right )^{2}
\leq \frac{\gamma_{2}}{\,\gamma_{1}}\,||u(x,\,0)||_{N}^{2}\quad \mbox{then}\quad
\frac{\gamma_{2}^{2}}{4\,\gamma_{1}^{2}}\,||u(x,\,0)||_{N}^{4} + 1 \leq 0
\end{eqnarray*}
which is a contradiction. Then $T_{N - 3}^{(n)}\not\rightarrow 0.$
By choosing $T_{N - 3}=T_{N - 3}(||u(x,\,0)||_{N}^{2})$
sufficiently small, and $T_{N - 3}$ not depending on $n,$ we
conclude that
\begin{eqnarray}
\label{e509}||v^{(n)}(\,\cdot\,,\,t)||_{\alpha}^{2}\leq c_{N -
3}^{2}\quad \mbox{for all} \quad 0\leq t\leq T_{N - 3}.
\end{eqnarray}
This shows that $T_{N - 3}^{(n)}\geq T_{N - 3}.$ Thus,
\begin{eqnarray*}
v\in L^{\infty}([0,\,T_{N - 3}]:\,H^{\alpha}(\mathbb{R}))\equiv
L^{\infty}([0,\,T_{N - 3}]:\,H^{N - 2}(\mathbb{R})).
\end{eqnarray*}
Now, denote by $0\leq T_{N - 3}^{*}\leq +\infty $ the maximal
number such that for all $0<t\leq T_{N - 3}^{*},$ $u=\wedge v\in
L^{\infty}([0,\,t]:\,H^{N}(\mathbb{R})).$ In particular, $T_{N -
3}\leq T_{N - 3}^{*}$ for all $N\geq 4.$ Thus, $T$ can be chosen
depending only on the norm of $u_{0}$ in $H^{3}(\mathbb{R}).$
Approximating $u_{0}$ by $\{u_{0}^{(j)}\}\in
C_{0}^{\infty}(\mathbb{R})$ such that $||u_{0} -
u_{0}^{(j)}||_{H^{N}(\mathbb{R})}\rightarrow 0$ as $j\rightarrow
+\infty.$ Let $u^{j}$ be a solution of \eqref{e204} with
$u^{(j)}(x,\,0)=u_{0}^{(j)}.$ According to the above argument,
there exists $T$ which is independent on $n$ but depending only on
$\sup_{j}||u_{0}^{(j)}||$ such that $u^{(j)}$ there exists on
$[0,\,T]$ and a subsequence $u^{(j)}\stackrel{j\rightarrow
+\infty}\longrightarrow u$
in $L^{\infty}([0,\,T]:\,H^{N}(\mathbb{R})).$\\
\\
As a consequence of Theorem 5.1 and 5.2 and its proof, one obtains the
following result.\\
\\
{\bf Corollary 5.3.} {\it Let $|\omega |<3\,\beta $ and let
$u_{0}\in H^{N}(\mathbb{R})$ with $N\geq 3$ such that
$u_{0}^{(j)}\rightarrow u_{0}$ in $H^{N}(\mathbb{R}).$ Let $u$ and
$u^{(j)}$ be the corresponding unique solutions given by Theorems
5.1 and 5.2 in $L^{\infty}([0,\,T]:\,H^{N}(\mathbb{R}))$ with $T$
depending only on $\sup_{j}||u_{0}^{(j)}||_{H^{3}(\mathbb{R})}$
such that}
\begin{eqnarray*}
u^{(j)}\stackrel{*}\rightharpoonup u\quad weakly \;on
\quad L^{\infty}([0,\,T]:\,H^{N}(\mathbb{R})),\\
u^{(j)}\rightarrow u\quad strongly \;on \quad L^{2}([0,\,T]:\,H^{N +
1}(\mathbb{R})).
\end{eqnarray*}
\renewcommand{\theequation}{\thesection.\arabic{equation}}
\setcounter{equation}{0}\section{Existence of Global Solutions}
Here, we will try to extend the local solution $u\in
L^{\infty}([0,\,T]:\,H^{N}(W_{0\;i\;0}))$ of \eqref{e204} obtained
in Theorem 5.2 to $t\geq 0.$ A standard way to obtain these
extensions consists into deducing global estimations for the
$H^{N}(W_{0\;i\;0})$-norm of $u$ in terms of the
$H^{N}(W_{0\;i\;0})$-norm of $u(x,\,0)=u_{0}(x).$ These
estimations are frequently based on conservation laws which
contain the $L^{2}$-norm of the solution and their spatial
derivatives. It is not possible to do the same to give a solution
of the problem of global existence because the difficulty here is
that the weight depends on the $x$ and $t$ variables. To solve our
problem we follow a different method using Leibniz's
rule like in the proof of Theorem 3.1 of Bona and Saut \cite{bo1}.\\
\\
{\bf Theorem 6.1.} {\it For $|\omega |<3\,\beta$ there exists a
global solution to \eqref{e204} in the space
$H^{s}(\mathbb{R})\cap H^{N}(W_{0\;i\;0})$ with
$N$ integer $\geq 3$ and $s\geq 2.$}\\
\\
{\it Proof.} The first part was proved in \cite{bo1}.
Differentiating \eqref{e204} $\alpha $-times (for $\alpha \geq 0$)
over $x\in \mathbb{R}$ leads to
\begin{eqnarray}
\label{e601}i\,u_{\alpha \,t} + i\,\beta \,u_{\alpha + 3} + \omega
\,u_{\alpha + 2} + (|u|^{2})_{\alpha }\,u + \sum_{m=1}^{\alpha -
1}{\alpha\choose m}(|u|^{2})_{\alpha - m}\,u_{m} +
|u|^{2}\,u_{\alpha } = 0.
\end{eqnarray}
Let $\xi = \xi(x,\,t),$ then multiplying \eqref{e601} by $\xi
\,\overline{u}_{\alpha } $ we have
\begin{eqnarray*}
\lefteqn{i\,\xi \,\overline{u}_{\alpha }\,u_{\alpha \,t} + i\,\beta
\,\xi \, \overline{u}_{\alpha }\,u_{\alpha + 3} + \omega \,\xi \,
\overline{u}_{\alpha }\,u_{\alpha + 2} + (|u|^{2})_{\alpha }\,
\xi \,u\,\overline{u}_{\alpha } } \\
&  & +\sum_{m=1}^{\alpha - 1}{\alpha\choose m} (|u|^{2})_{\alpha -
m}\,\xi \,u_{m}\,\overline{u}_{\alpha} + \xi \,|u|^{2}\,|u_{\alpha
}|^{2}  =  0
\end{eqnarray*}
and
\begin{eqnarray*}
\lefteqn{-\,i\,\xi \,u_{\alpha }\,\overline{u}_{\alpha \,t} -
i\,\beta \,\xi \,u_{\alpha }\,\overline{u}_{\alpha + 3} + \omega
\,\xi \,u_{\alpha }\,\overline{u}_{\alpha + 2} + (|u|^{2})_{\alpha
}\,
\xi \,\overline{u}\,u_{\alpha } } \\
&  & +\sum_{m=1}^{\alpha - 1}{\alpha\choose m}(|u|^{2})_{\alpha -
m}\,\xi \,\overline{u}_{m}\,u_{\alpha} + \xi \,|u|^{2}\,|u_{\alpha
}|^{2}  =  0. \qquad (\mbox{applying conjugate})
\end{eqnarray*}
Subtracting and integrating over $x\in \mathbb{R}$ we have
\begin{eqnarray}
\label{e602}\lefteqn{i\,\partial_{t}\int_{\mathbb{R}}\xi
\,|u_{\alpha }|^{2}dx + i\,\beta \int_{\mathbb{R}}\xi
\,\overline{u}_{\alpha }\,u_{\alpha + 3}dx + i\,\beta
\int_{\mathbb{R}}\xi \,u_{\alpha}\,\overline{u}_{\alpha + 3}dx +
\omega \int_{\mathbb{R}}\xi \,\overline{u}_{\alpha}\,u_{\alpha +
2}dx } \\
&  & -\;\omega \int_{\mathbb{R}}\xi
\,u_{\alpha}\,\overline{u}_{\alpha + 2}dx +
2\,i\,Im\int_{\mathbb{R}}\xi \,(|u|^{2})_{\alpha
}\,u\,\overline{u}_{\alpha }dx + 2\,i\sum_{m=1}^{\alpha -
1}{\alpha\choose m}Im\int_{\mathbb{R}}\xi \,(|u|^{2})_{\alpha - m}
\,u_{m}\,\overline{u}_{\alpha}dx = 0. \nonumber
\end{eqnarray}
Each term is calculated separately, integrating by parts in the
second term we have
\begin{eqnarray*}
\int_{\mathbb{R}}\xi\,\overline{u}_{\alpha}\,u_{\alpha + 3}dx =
\int_{\mathbb{R}}\partial^{2}\xi
\,\overline{u}_{\alpha}\,u_{\alpha + 1}dx +
2\int_{\mathbb{R}}\partial\xi\,|u_{\alpha + 1}|^{2}dx +
\int_{\mathbb{R}}\xi\,\overline{u}_{\alpha + 2}\,u_{\alpha + 1}dx.
\end{eqnarray*}
The other terms are calculated in a similar way. Hence in
\eqref{e602}
\begin{eqnarray*}
\lefteqn{\partial_{t}\int_{\mathbb{R}}\xi \,|u_{\alpha}|^{2}dx -
\beta \int_{\mathbb{R}}\partial^{3}\xi\,|u_{\alpha}|^{2}dx +
3\,\beta \int_{\mathbb{R}}\partial \xi\,|u_{\alpha + 1}|^{2}dx -
2\,\omega \,Im\int_{\mathbb{R}}\partial
\xi\,\overline{u}_{\alpha}\,u_{\alpha + 1}dx } \\
&  & -\int_{\mathbb{R}}\partial_{t}\xi\,|u_{\alpha}|^{2}dx +
2\,Im\int_{\mathbb{R}}\xi \,(|u|^{2})_{\alpha
}\,u\,\overline{u}_{\alpha }dx + 2\sum_{m=1}^{\alpha -
1}{\alpha\choose m}Im\int_{\mathbb{R}}\xi \,(|u|^{2})_{\alpha - m}
\,u_{m}\,\overline{u}_{\alpha}dx = 0
\end{eqnarray*}
such that
\begin{eqnarray*}
\lefteqn{\partial_{t}\int_{\mathbb{R}}\xi \,|u_{\alpha}|^{2}dx -
\beta \int_{\mathbb{R}}\partial^{3}\xi\,|u_{\alpha}|^{2}dx +
3\,\beta \int_{\mathbb{R}}\partial \xi\,|u_{\alpha + 1}|^{2}dx +
2\,Im\int_{\mathbb{R}}(|u|^{2})_{\alpha}\,\xi
\,u\,\overline{u}_{\alpha}dx}\\
&  & - \int_{\mathbb{R}}\partial_{t}\xi\,|u_{\alpha}|^{2}dx +
2\sum_{m=1}^{\alpha - 1}{\alpha\choose m} Im\int_{\mathbb{R}}\xi
\,(|u|^{2})_{\alpha - m}
\,u_{m}\,\overline{u}_{\alpha}dx \\
& = & 2\,\alpha \,Im\int_{\mathbb{R}}\partial \xi\,
\overline{u}_{\alpha}\,u_{\alpha + 1}dx \leq |\omega
|\int_{\mathbb{R}}\partial \xi\,|u_{\alpha}|^{2}dx + |\omega
|\int_{\mathbb{R}}\partial \xi\,|u_{\alpha + 1}|^{2}dx.
\end{eqnarray*}
Hence
\begin{eqnarray}
\lefteqn{\partial_{t}\int_{\mathbb{R}}\xi \,|u_{\alpha}|^{2}dx +
\int_{\mathbb{R}}[\,3\,\beta - |\omega |\,]\,\partial
\xi\;|u_{\alpha + 1}|^{2}dx - \int_{\mathbb{R}}[\,\partial_{t}\xi +
\beta \,\partial^{3}\xi
+ |\omega|\,\partial \xi\,]\,|u_{\alpha}|^{2}dx}\nonumber \\
\label{e603}&  & +\;2\,Im\int_{\mathbb{R}}(|u|^{2})_{\alpha}\,\xi
\,u\,\overline{u}_{\alpha}dx + 2\sum_{m=1}^{\alpha -
1}{\alpha\choose m}Im\int_{\mathbb{R}}\xi \,(|u|^{2})_{\alpha - m}
u_{m}\,\overline{u}_{\alpha}dx\leq 0.
\end{eqnarray}
But
\begin{eqnarray*}
(|u|^{2})_{\alpha} & = & (\,u\,\overline{u}\,)_{\alpha}  =
\sum_{k=0}^{\alpha}{\alpha\choose k}u_{\alpha - k}\,\overline{u}_{k}
=  \overline{u}\,u_{\alpha} + \sum_{k=1}^{\alpha - 1}{\alpha\choose
k}u_{\alpha - k} \,\overline{u}_{k} + u\,\overline{u}_{\alpha}
\end{eqnarray*}
then
\begin{eqnarray*}
(|u|^{2})_{\alpha}\,u\,\overline{u}_{\alpha} =
|u|^{2}\,|u_{\alpha}|^{2} + \sum_{k=1}^{\alpha - 1}{\alpha\choose
k}u_{\alpha - k}\, \overline{u}_{k}\,u\,\overline{u}_{\alpha} +
u^{2}\,\overline{u}_{\alpha}^{2}
\end{eqnarray*}
hence
\begin{eqnarray}
\lefteqn{2\;Im\int_{\mathbb{R}}(\,|u|^{2}\,)_{\alpha}\,\xi
\,u\,\overline{u}_{\alpha}dx = 2\sum_{k=1}^{\alpha -
1}{\alpha\choose k}Im\int_{\mathbb{R}}\xi \,u_{\alpha -
k}\,\overline{u}_{k}\,u\,\overline{u}_{\alpha}dx
+ 2\;Im\int_{\mathbb{R}}\xi \,u^{2}\,\overline{u}_{\alpha}^{2}dx } \nonumber \\
& \leq & 2\sum_{k=1}^{\alpha - 1}{\alpha\choose
k}\int_{\mathbb{R}}\xi \,|u_{\alpha -
k}|\,|u_{k}|\,|u|\,|u_{\alpha}|dx + 2\int_{\mathbb{R}}\xi\,|u|^{2}
\,|u_{\alpha}|^{2}dx \nonumber \\
& \leq & 2\sum_{k=1}^{\alpha - 1}{\alpha\choose
k}\int_{\mathbb{R}}\xi \,|u_{\alpha -
k}|\,|u_{k}|\,|u|\,|u_{\alpha}|dx +
2\,||u||_{L^{\infty}(\mathbb{R})}^{2}\int_{\mathbb{R}}\xi
\,|u_{\alpha}|^{2}dx
\nonumber \\
\label{e604}& \leq &
2\,||u||_{L^{\infty}(\mathbb{R})}\sum_{k=1}^{\alpha - 1}
{\alpha\choose k}\int_{\mathbb{R}}\xi \,|u_{\alpha -
k}|\,|u_{k}|\,|u_{\alpha}|dx +
2\,||u||_{L^{\infty}(\mathbb{R})}^{2}\int_{\mathbb{R}}\xi
\,|u_{\alpha}|^{2}dx
\end{eqnarray}
hence in \eqref{e603} we have
\begin{eqnarray*}
\lefteqn{\partial_{t}\int_{\mathbb{R}}\xi \,|u_{\alpha}|^{2}dx +
\int_{\mathbb{R}}[3\,\beta - |\omega |\,]\,\partial \xi\,|u_{\alpha
+ 1}|^{2}dx \leq \int_{\mathbb{R}}[\partial_{t}\xi + \beta
\,\partial^{3}\xi + |\omega|\,\partial \xi
+ c\,\xi \,]\,|u_{\alpha}|^{2}dx }\\
&  & +\,2\,c\sum_{k=1}^{\alpha - 1}{\alpha\choose
k}\int_{\mathbb{R}}\xi \,|u_{\alpha -
k}|\,|u_{k}|\,|u|\,|u_{\alpha}|dx - 2\sum_{m=1}^{\alpha -
1}{\alpha\choose m}Im\int_{\mathbb{R}}\xi \,(|u|^{2})_{\alpha - m}
\,u_{m}\,\overline{u}_{\alpha}dx.
\end{eqnarray*}
Using \eqref{e203}, Gagliardo-Nirenberg's inequality and standard
estimates we get
\begin{eqnarray}
\label{e605}\partial_{t}\int_{\mathbb{R}}\xi\,|u_{\alpha}|^{2}dx +
[3\,\beta - |\omega |\,]\int_{\mathbb{R}}\partial \xi\,|u_{\alpha
+ 1}|^{2}\,dx \leq c\int_{\mathbb{R}}\xi \,|u_{\alpha}|^{2}dx.
\end{eqnarray}
Integrating \eqref{e605} in $t\in [0,\,T_{max}=T]$ we obtain
\begin{eqnarray*}
\int_{\mathbb{R}}\xi\,|u_{\alpha}|^{2}dx + [3\,\beta - |\omega
|\,]\int_{0}^{t}\int_{\mathbb{R}}\partial \xi\,|u_{\alpha +
1}|^{2}dx\,ds \leq ||u_{0}(x)||_{\alpha}^{2} +
\int_{0}^{t}\left(c\int_{\mathbb{R}}\xi
\,|u_{\alpha}|^{2}dx\right)ds,
\end{eqnarray*}
where
\begin{eqnarray*}
\int_{\mathbb{R}}\xi \,|u_{\alpha}|^{2}dx  \leq
||u_{0}(x)||_{\alpha}^{2} + \int_{0}^{t}\left(c\int_{\mathbb{R}}\xi
\,|u_{\alpha}|^{2}dx\right)ds.
\end{eqnarray*}
Using Gronwall's inequality
\begin{eqnarray*}
\int_{\mathbb{R}}\xi \,|u_{\alpha}|^{2}dx \leq
||u_{0}(x)||_{\alpha}^{2}\,e^{c\,t}\leq
||u_{0}(x)||_{\alpha}^{2}\,e^{c\,T}
\end{eqnarray*}
it follows that
\begin{eqnarray*}
\int_{\mathbb{R}}\xi \,|u_{\alpha}|^{2}dx \leq
c=c(T,\,||u_{0}(x)||_{\alpha}^{2}).
\end{eqnarray*}
Then for any $T=T_{max}>0$ there exists
$c=c(T,\,||u_{0}(x)||_{\alpha}^{2})$ such that
\begin{eqnarray*}
||u||_{\alpha}^{2} + [3\,\beta - |\omega
|\,]\int_{0}^{t}\int_{\mathbb{R}}\partial\xi\,|u_{\alpha +
1}|^{2}dx\,ds \leq c.
\end{eqnarray*}
This concludes the proof.
\renewcommand{\theequation}{\thesection.\arabic{equation}}
\setcounter{equation}{0}\section{Persistence Theorem} As a
starting point for the a priori gain of regularity results that
will be discussed in the next section, we need to develop some
estimates for solutions of the equation \eqref{e204} in weighted
Sobolev norms. The existence of these weighted estimates is often
called the persistence of a property of the initial data $u_{0}.$
We show that if $u_{0}\in H^{3}(\mathbb{R})\cap
H^{L}(W_{0\;i\;0})$ for $L\geq 0,$ $i\geq 1,$ then the solution
$u(\,\cdot \,,\,t)$ evolves in $H^{L}(W_{0\;i\;0})$ for $t\in
[0,\,T].$ The time interval of that persistence is at least as
long as the interval guaranteed
by the existence Theorem 5.2.\\
\\
{\bf Theorem 7.1} (Persistence). {\it Let $|\omega |<3\,\beta$ and
let $i\geq 1$ and $L\geq 0$ be non-negative integers, $0<T<+\infty.$
Assume that $u$ is the solution to \eqref{e204} in
$L^{\infty}([0,\,T]:\,H^{3}(\mathbb{R}))$ with initial data
$u_{0}(x)=u(x,\,0)\in H^{3}(\mathbb{R}).$ If $u_{0}(x) \in
H^{L}(W_{0\;i\;0})$ then}
\begin{eqnarray}
\label{e701}u\in L^{\infty}\left([0,\,T]:\,H^{3}(\mathbb{R})
\cap H^{L}(W_{0\;i\;0})\right)\\
\label{e702}\int_{0}^{T}\int_{\mathbb{R}}|\partial^{L +
1}u(x,\,t)|^{2}\,\eta\,dx\,dt<+\infty
\end{eqnarray}
{\it where $\sigma$ is arbitrary, $\eta\in W_{\sigma\;i\;0}$ for $i\geq 1.$}\\
\\
{\it Proof.} We use induction on $\alpha .$ Let
\begin{eqnarray*}
u\in L^{\infty}\left([0,\,T]:\,H^{3}(\mathbb{R})\cap
H^{\alpha}(W_{0\;i\;0})\right) \quad \mbox{for}\quad 0\leq \alpha
\leq L.
\end{eqnarray*}
We derive formally some a priori estimate for the solution where the
bound, involves only the norms of $u$ in
$L^{\infty}([0,\,T]:\,H^{3}(\mathbb{R}))$ and the norms of $u_{0}$
in $H^{3}(W_{0\;i\;0}).$ We do this by approximating $u(x,\,t)$
through smooth solutions and the weight functions by smooth bounded
functions. By Theorem 5.2, we have
\begin{eqnarray*}
u(x,\,t)\in L^{\infty}([0,\,T]:\,H^{N}(\mathbb{R}))\quad
\mbox{with}\quad N=\mbox{max}\{L,\,3\}.
\end{eqnarray*}
In particular, $u_{j}(x,\,t)\in L^{\infty}([0,\,T]\times
\mathbb{R})$ for $0\leq j\leq N - 1.$ To obtain \eqref{e701} and
\eqref{e702} there are two ways of approximation. We approximate
general solutions by smooth solutions, and we approximate general
weight functions by bounded weight functions. The first of these
procedure has already been discussed,
so we shall concentrate on the second.\\
Given a smooth weight function $\eta (x)\in W_{\sigma ,\;i -
1,\;0}$ with $\sigma >0,$ we take a sequence $\eta^{\nu}(x)$ of
smooth bounded weight functions approximating $\eta (x)$ from
below, uniformly on any half line $(-\infty,\,c).$ Define the
weight functions for the $\alpha$-th induction step as
\begin{eqnarray*}
\xi _{\nu}=\frac{1}{(3\,\beta -
|\omega|)}\int_{-\infty}^{x}\eta^{\nu}(y,\,t)\,dy
\end{eqnarray*}
then the $\xi_{\nu}$ are bounded weight functions which approximate
a desired weight function $\xi\in W_{0\;i\;0}$ from below, uniformly
on a compact set.  For $\alpha =0,$ multiplying \eqref{e204} by
$\xi_{\nu}\,\overline{u},$ we have
\begin{eqnarray*}
&  & i\,\xi_{\nu}\,\overline{u}\,u_{t} + i\,\beta\,\xi_{\nu}\,
\overline{u}\,u_{3} + \omega \,\xi_{\nu}\,\overline{u}\,u_{2} +
\xi_{\nu}\,|u|^{4}=0 \\
&  & -\,i\,\xi_{\nu}\,u\,\overline{u}_{t} - i\,\beta \,\xi
_{\nu}\,u\,\overline{u}_{3} + \omega\,\xi
_{\nu}\,u\,\overline{u}_{2} + \xi_{\nu}\,|u|^{4}=0. \quad
\mbox{(applying conjugate)}
\end{eqnarray*}
Subtracting and integrating over $x\in \mathbb{R}$ we have
\begin{eqnarray}
&  & i\,\partial_{t}\int_{\mathbb{R}}\xi_{\nu}\,|u|^{2}dx -
i\int_{\mathbb{R}}\partial_{t}\xi_{\nu}\,|u|^{2}dx + i\,\beta
\int_{\mathbb{R}}\xi_{\nu}\,\overline{u}\,u_{3}dx
+ i\,\beta \int_{\mathbb{R}}\xi_{\nu}\,u\,\overline{u}_{3}dx \nonumber \\
\label{e703}&  & +\,\omega
\int_{\mathbb{R}}\xi_{\nu}\,\overline{u}\,u_{2}dx - \omega
\int_{\mathbb{R}}\xi_{\nu}\,u\,\overline{u}_{2}dx = 0.
\end{eqnarray}
Each term is treated separately, integrating by parts in the third
term we have
\begin{eqnarray*}
\int_{\mathbb{R}}\xi_{\nu}\,\overline{u}\,u_{3}dx & = &
\int_{\mathbb{R}}\partial^{2}\xi_{\nu}\,\overline{u}\,u_{1}dx +
2\int_{\mathbb{R}}\partial\xi_{\nu}\,|u_{1}|^{2}dx +
\int_{\mathbb{R}}\xi_{\nu}\,\overline{u}_{2}\,u_{1}dx.
\end{eqnarray*}
The other terms are calculated in a similar way. Hence in
\eqref{e703} we have
\begin{eqnarray*}
\lefteqn{\partial_{t}\int_{\mathbb{R}}\xi_{\nu}\,|u|^{2}dx -
\int_{\mathbb{R}}\partial_{t}\xi_{\nu}\,|u|^{2}dx - \beta
\int_{\mathbb{R}}\partial^{3}\xi_{\nu}\,|u|^{2}dx + 3\,\beta
\int_{\mathbb{R}}\partial \xi_{\nu}\,|u_{1}|^{2}dx } \\
& = & 2\,\omega \,Im\int_{\mathbb{R}}\partial
\xi_{\nu}\,\overline{u}\,u_{1}dx \leq  |\omega
|\int_{\mathbb{R}}\partial \xi_{\nu}\,|u|^{2}dx + |\omega
|\int_{\mathbb{R}}\partial \xi_{\nu}\,|u_{1}|^{2}dx.
\end{eqnarray*}
Then, using \eqref{e203} we obtain
\begin{eqnarray*}
\lefteqn{\partial_{t}\int_{\mathbb{R}}\xi_{\nu}\,|u|^{2}dx +
\int_{\mathbb{R}}[3\,\beta - |\omega |]\,\partial
\xi_{\nu}\,|u_{1}|^{2}dx }\\
& \leq &\int_{\mathbb{R}}[\partial_{t}\xi_{\nu} + \beta
\,\partial^{3}\xi_{\nu} + |\omega |\,\partial \xi_{\nu}]\,|u|^{2}dx
\leq c\int_{\mathbb{R}}\xi_{\nu}\,|u|^{2}dx
\end{eqnarray*}
thus
\begin{eqnarray*}
\partial_{t}\int_{\mathbb{R}}\xi_{\nu}\,|u|^{2}dx \leq
c\int_{\mathbb{R}}\xi_{\nu}\,|u|^{2}dx.
\end{eqnarray*}
We apply Gronwall's Lemma to conclude that
\begin{eqnarray}
\label{e704}\partial_{t}\int_{\mathbb{R}}\xi_{\nu}\,|u|^{2}dx \leq
c(T,\,||u_{0}||).
\end{eqnarray}
for $0\leq t\leq T,$ and $c$ not depending on $\beta >0,$
the weighted estimate remains true for $\beta \rightarrow 0.$\\
Now, we assume that the result is true for $(\alpha - 1)$ and we
prove that it is true for $\alpha .$ To prove this, we start from
the main inequality \eqref{e301} with $\xi$ and $\eta$ given by
$\xi_{\nu}$ and $\eta_{\nu}$ respectively.
\begin{eqnarray*}
&  & \partial_{t}\int_{\mathbb{R}}\xi _{\nu}\,|u_{\alpha}|^{2}dx +
\int_{\mathbb{R}}\eta_{\nu}\,|u_{\alpha + 1}|^{2}dx +
\int_{\mathbb{R}}\theta_{\nu}\,|u_{\alpha}|^{2}dx +
\int_{\mathbb{R}}R_{\alpha}dx\leq 0
\end{eqnarray*}
where
\begin{eqnarray*}
\eta_{\nu} & = & (3\beta - |\omega |\,)\,\partial \xi_{\nu}
\qquad \mbox{for}\qquad |\omega |<3\;\beta \\
\theta_{\nu} & = & -\;[\,\partial_{t}\xi_{\nu} + \beta
\,\partial^{3}\xi_{\nu} + |\omega|\,\partial \xi_{\nu} +
c_{0}\,\xi_{\nu} \,]\qquad \mbox{where}\quad
c_{0}=||u||_{L^{\infty}(\mathbb{R})}^{2}\\
R_{\alpha} & = & R_{\alpha}(|u_{\alpha }|,\,|u_{\alpha -
1}|,\,\ldots \,)
\end{eqnarray*}
then
\begin{eqnarray*}
\lefteqn{\partial_{t}\int_{\mathbb{R}}\xi _{\nu}\,|u_{\alpha}|^{2}dx
+ \int_{\mathbb{R}}\eta_{\nu} \,|u_{\alpha + 1}|^{2}dx \leq
-\int_{\mathbb{R}}\theta_{\nu}\,|u_{\alpha}|^{2}dx -
\int_{\mathbb{R}}R_{\alpha}dx }\\
& \leq & \left|-\int_{\mathbb{R}}\theta_{\nu}\,|u_{\alpha}|^{2}dx
- \int_{\mathbb{R}}R_{\alpha}dx\,\right|\leq
\int_{\mathbb{R}}|\theta_{\nu}|\,|u_{\alpha}|^{2}dx +
\int_{\mathbb{R}}|R_{\alpha}|dx.
\end{eqnarray*}
Using \eqref{e203} in the first part of the right hand side we
obtain
\begin{eqnarray*}
\int_{\mathbb{R}}\theta_{\nu}\,|u_{\alpha}|^{2}dx \leq
c\int_{\mathbb{R}}\xi_{\nu}\,|u_{\alpha}|^{2}dx
\end{eqnarray*}
thus
\begin{eqnarray}
\label{e705}\partial_{t}\int_{\mathbb{R}}\xi
_{\nu}\,|u_{\alpha}|^{2}dx + \int_{\mathbb{R}}\eta_{\nu}
\,|u_{\alpha + 1}|^{2}dx \leq
c\int_{\mathbb{R}}\xi_{\nu}\,|u_{\alpha}|^{2}dx +
\int_{\mathbb{R}}|R_{\alpha}|dx.
\end{eqnarray}
According to \eqref{e308}, $\int_{\mathbb{R}}R_{\alpha}\,dx$
contains a term of the form
\begin{eqnarray}
\label{e706}
\int_{\mathbb{R}}\xi_{\nu}\,u_{\nu_{1}}\,\overline{u}_{\nu_{2}}\,
\overline{u}_{\alpha}dx.
\end{eqnarray}
We estimate the term
\begin{eqnarray}
\label{e707}\int_{\mathbb{R}}\xi_{\nu}\,u_{\nu_{1}}\,\overline{u}_{\nu_{2}}\,
\overline{u}_{\alpha}\,dx\quad \mbox{for}\quad \nu_{1} +
\nu_{2}=\alpha .
\end{eqnarray}
Let $\nu_{2}\leq \alpha - 2.$ Integrating by parts one time in
\eqref{e707} we have
\begin{eqnarray*}
\int_{\mathbb{R}}\xi_{\nu}\,u_{\nu_{1}}\,\overline{u}_{\nu_{2}}\,
\overline{u}_{\alpha}\,dx & = &
-\int_{\mathbb{R}}\partial\xi_{\nu}\,u_{\nu_{1}}\,\overline{u}_{\nu_{2}}\,
\overline{u}_{\alpha - 1}\,dx -
\int_{\mathbb{R}}\xi_{\nu}\,u_{\nu_{1} +
1}\,\overline{u}_{\nu_{2}}\,
\overline{u}_{\alpha - 1}\,dx\\
&  & -\int_{\mathbb{R}}\xi_{\nu}\,u_{\nu_{1}}\,\overline{u}_{\nu_{2}
+ 1}\, \overline{u}_{\alpha - 1}\,dx.
\end{eqnarray*}
We estimates the first term in the right hand side in \eqref{e707}.
Using Holder's inequality and standard estimates we obtain
\begin{eqnarray}
\label{e709}c\,\left[\left(\int_{\mathbb{R}}\xi_{\nu}\,|u_{\nu_{2} +
1}|^{2}dx\right)^{1/2} +
\left(\int_{\mathbb{R}}\xi_{\nu}\,|u_{\nu_{2}}|^{2}dx\right)^{1/2}
\right]\left(\int_{\mathbb{R}}\xi_{\nu}\,|u_{\alpha -
1}|^{2}dx\right)^{1/2}
\end{eqnarray}
where \eqref{e709} is bounded by hypothesis. The other terms are
estimates in a similar way. Now suppose that $\nu_{1}=\nu_{2}=
\alpha - 1,$ then in \eqref{e707} we have
\begin{eqnarray*}
\int_{\mathbb{R}}\xi_{\nu}\,u_{\alpha - 1}\,\overline{u}_{\alpha -
1}\,\overline{u}_{\alpha}dx,
\end{eqnarray*}
hence
\begin{eqnarray*}
\left |\int_{\mathbb{R}}\xi_{\nu}\,|u_{\alpha -
1}|^{2}\,\overline{u}_{\alpha}dx \right |\leq ||u_{\alpha -
1}||_{L^{\infty}(\mathbb{R})}\,
\left(\int_{\mathbb{R}}\xi_{\nu}\,|u_{\alpha -
1}|^{2}dx\right)^{1/2}\,
\left(\int_{\mathbb{R}}\xi_{\nu}\,|u_{\alpha}|^{2}dx\right)^{1/2}
\end{eqnarray*}
where $||u_{\alpha - 1}||_{L^{\infty}(\mathbb{R})}$ is bounded by
hypothesis, and the estimate is complete. In a similar way we
estimate all the other terms of $R_{\alpha}.$ Using these estimates
in \eqref{e705} and applying Gronwall's argument, we obtain for
$0\leq t\leq T$
\begin{eqnarray*}
\partial_{t}\int_{\mathbb{R}}\xi_{\nu}\,|u_{\alpha}|^{2}dx +
\int_{\mathbb{R}}\eta_{\nu}\,|u_{\alpha + 1}|^{2}dx \leq
c_{0}\,e^{c_{1}\,t}\,\left(\int_{\mathbb{R}}\xi_{\nu}\,
|\partial^{\alpha}u_{0}(x)|^{2}dx + 1\right)
\end{eqnarray*}
where $c_{0}$ and $c_{1}$ are independent of $\nu$ and such that
letting the parameter $\nu \rightarrow 0$ the desired estimate
\eqref{e702} is obtained.
\renewcommand{\theequation}{\thesection.\arabic{equation}}
\setcounter{equation}{0}\section{Main Theorem} In this section we
state and prove our main theorem, which states that if the initial
data $u(x,\,0)$ decays faster than polynomially on
$\mathbb{R}^{+}=\{x\in \mathbb{R}:\;x>0\}$ and possesses certain
initial Sobolev regularity, then the solution $u(x,\,t)\in
C^{\infty}$ for all $t>0.$ \\
\\
If $\eta$ is an arbitrary weight function in $W_{\sigma\;i\;k},$
then by Lemma 3.2, there exists $\xi\in W_{\sigma,\;i + 1,\;k}$
which satisfies \eqref{e301}. For the main theorem, we take $4\leq
\alpha \leq L + 2.$ For $\alpha \leq L + 4,$ we take
\begin{eqnarray}
\label{e801}\eta \in W_{\sigma ,\,L - \alpha - 2,\,\alpha - 3}\;
\Longrightarrow \;\xi \in W_{\sigma ,\,L - \alpha - 3,\,\alpha - 3}.
\end{eqnarray}
{\bf Lemma 8.1}(Estimate of error terms). {\it Let $4\leq \alpha
\leq L + 2$ and the weight functions be chosen as in \eqref{e801},
then}
\begin{eqnarray}
\label{e802}\left|\int_{0}^{T}\int_{\mathbb{R}}(\theta
\,|u_{\alpha}|^{2} + R_{\alpha})dx\,dt\right| \leq c,
\end{eqnarray}
{\it where $c$ depends only on the norms of $u$ in }
\begin{eqnarray*}
L^{\infty}([0,\,T]:\,H^{\beta}(W_{\sigma,\,L - \beta + 3,\,\beta - 3}))
\cap L^{2}([0,\,T]:\,H^{\beta + 1}(W_{\sigma,\,L - \beta + 2,\,\beta - 3}))
\end{eqnarray*}
{\it for $3\leq \beta \leq \alpha - 1,$ and the norms of
$u$ in} $L^{\infty}([0,\,T]:\,H^{3}(W_{0\;L\;0})).$\\
\\
{\it Proof.} We must estimate both $R_{\alpha}$ and $\theta .$ We
begin with a term in $R_{\alpha}$ of the form
\begin{eqnarray}
\label{e803}\xi \,|u_{\nu_{1}}|\,\,|u_{\nu_{2}}|\,\,|u_{\alpha}|
\end{eqnarray}
assuming that $\nu_{1}\leq \alpha - 2.$ \\
\\
By the induction
hypothesis, $u$ is bounded in
$L^{\infty}([0,\,T]:\,H^{\beta}(W_{\sigma,\,L - (\beta - 3)^{+},
\,(\beta - 3)^{+}}))$ for $0\leq \beta \leq \alpha - 1.$ By Lemma
2.1,
\begin{eqnarray}
\label{e804}\sup_{t\geq 0} \,\sup_{x\in\mathbb{R}}\,\zeta
\,|u_{\beta}|^{2}<+\infty
\end{eqnarray}
for $0\leq \beta \leq \alpha - 2$ and $\zeta\in W_{\sigma ,\,L -
(\beta - 2)^{+},\,(\beta - 2)^{+}}.$ We estimate $|u_{\nu_{1}}|$
using \eqref{e804}. We estimate $|u_{\nu_{2}}|$ and $|u_{\alpha}|$
using the weighted $L^{2}$ bounds
\begin{eqnarray}
\label{e805}\int_{0}^{T}\int_{\mathbb{R}}\zeta
\,|u_{\nu_{2}}|^{2}dx\,dt<+\infty \quad \mbox{for}\quad \zeta \in
W_{\sigma ,\,L - (\nu_{2} - 3)^{+},\,(\nu_{2} - 4)^{+}}
\end{eqnarray}
and the same with $\nu_{2}$ replaced by $\alpha .$ It suffices to check the
powers to
$t,$ the powers of $x$ as $x\rightarrow +\infty $ and the exponential of
$x$ as $x\rightarrow -\infty .$\\
\\
For $x>1.$ In the \eqref{e803} term, the factor $\xi $ contributed
according to \eqref{e801}
\begin{eqnarray*}
\xi(x,\,t) = t^{\alpha - 3}\,x^{(L - \alpha + 3)}\,t^{-(\alpha -
3)}\,x^{-(L - \alpha + 3)}\xi(x,\,t)\leq c_{2}\,t^{\alpha -
3}\,x^{(L - \alpha + 3)}\quad (\mbox{using} \eqref{e203})
\end{eqnarray*}
then $\xi \,|u_{\nu_{1}}|\,|u_{\nu_{2}}|\,|u_{\alpha}|\leq
c_{2}\,t^{\alpha - 3}\,x^{(L - \alpha +
3)}|u_{\nu_{1}}|\,|u_{\nu_{2}}|\,|u_{\alpha}|.$ Moreover
\begin{eqnarray*}
|u_{\nu_{1}}|\,|u_{\nu_{2}}|\,|u_{\alpha}| & = &
t^{\frac{(\nu_{1} - 2)^{+}}{2}}\,
x^{\frac{L - (\nu_{1} - 2)^{+}}{2}}\,
t^{\frac{-(\nu_{1} - 2)^{+}}{2}}\,
x^{\frac{(L - (\nu_{1} - 2)^{+})}{2}}\,|u_{\nu_{1}}|\times \\
&  & t^{\frac{(\nu_{2} - 4)^{+}}{2}}\, x^{\frac{L - (\nu_{2} -
3)^{+}}{2}}\,t^{\frac{-(\nu_{2} - 4)^{+}}{2}}\,x^{\frac{(L -
(\nu_{2} - 3)^{+})}{2}}\,
|u_{\nu_{2}}|\times \\
& & t^{\frac{(\alpha - 4)^{+}}{2}}\, x^{\frac{L - (\alpha -
3)^{+}}{2}}\,t^{\frac{-(\alpha - 4)^{+}}{2}}\,x^{\frac{(L -
(\alpha - 3)^{+})}{2}}\,|u_{\alpha}|.
\end{eqnarray*}
tt follows that
\begin{eqnarray}
\lefteqn{\xi
\,|u_{\nu_{1}}|\,|u_{\nu_{2}}|\,|u_{\alpha}|}\nonumber\\
\label{e806}&  &  \leq c_{2}\,t^{M}\,x^{T}\,t^{\frac{(\nu_{1} -
2)^{+}}{2}}\, x^{\frac{L - (\nu_{1} - 2)^{+}}{2}}\,|u_{\nu_{1}}|\,
t^{\frac{(\nu_{2} - 4)^{+}}{2}}\, x^{\frac{L - (\nu_{2} -
3)^{+}}{2}}\,|u_{\nu_{2}}|\,t^{\frac{(\alpha - 4)^{+}}{2}}\,
x^{\frac{L - (\alpha - 3)^{+}}{2}}\,\,|u_{\alpha}|\quad
\end{eqnarray}
where
\begin{eqnarray*}
M=\alpha - 3 - \frac{1}{2}(\nu_{1} - 2)^{+} -
\frac{1}{2}(\nu_{2} - 4)^{+} - \frac{1}{2}(\alpha - 4)^{+}
\end{eqnarray*}
and
\begin{eqnarray*}
T=(L - \alpha + 3) - \frac{1}{2}(L - (\alpha - 3)^{+}) -
\frac{1}{2}(L - (\nu_{2} - 3)^{+}) - \frac{1}{2}(L - \nu_{1} - 2)^{+}).
\end{eqnarray*}
{\it Claim.} $M\geq 0$ is large enough, that the extra power of $t$
can be omitted
\begin{eqnarray*}
2\,M & = & 2\,\alpha - 6 - (\nu_{1} - 2)^{+} - (\nu_{2} - 4)^{+} - (\alpha - 4)^{+}\\
& = & \alpha - 2 - (\nu_{1} - 2)^{+} - (\nu_{2} - 4)^{+}\\
& = & \alpha - 2 - \nu_{1} + 2 - \nu_{2} + 4
 =  \alpha + 4 - (\nu_{1} + \nu_{2})\\
& = & \alpha + 4 - \alpha   =  4\geq 0.
\end{eqnarray*}
{\it Claim.} $T\leq 0$ is such that the extra power of $t$ can be
omitted.
\begin{eqnarray*}
2\,T & = & 2\,L - 2\,\alpha + 6 - L + (\alpha - 3)^{+} - L +
(\nu_{2} - 3)^{+} - L + (\nu_{1} - 2)^{+}\\
& = & -\,L - \alpha + \nu_{1} + \nu_{2} - 2
 =  -\,L - \alpha + \alpha - 2\\
& = & -(L + 2) \leq 0.
\end{eqnarray*}
Now, we study the behavior as $x\rightarrow -\infty .$ Since each
factor $u_{\nu_{j}}$($j=1,\,2$) must grow slower that an
exponential $e^{\sigma'\,|x|}$ and $\xi$ decays as an exponential
$e^{-\sigma \,|x|},$ we simply need to choose the appropriate
relationship $\sigma$ and $\sigma'$ at each induction step. The
analysis will be completed with the case where $\nu_{1}\geq \alpha
- 1.$ Then, in \eqref{e309}, if $2(\alpha - 1)\leq \alpha,$ but
$\alpha \geq 3.$ So this possibility is impossible. For $x<1$ the
estimate is similar, except for an exponential weight. The
analysis of all terms of $R_{\alpha}$ is estimated
in a similar form. This completes the estimate of $R_{\alpha}.$ \\
Now, we estimate the term $\theta \,|u_{\alpha}|^{2}$ where
$\theta $ is given in \eqref{e301}. We have that $\theta $
involves derivatives of $u$ only up to order one, and hence,
$\theta \,|u_{\alpha}|^{2}$ is a sum of terms of the same type
which we have already encountered in $R_{\alpha}.$ So, its
integral can be bounded in the same type. Indeed, \eqref{e301}
shows that $\theta$ depends on $\xi_{t},$ $\partial^{3}\xi$ and
derivatives of lower order. By using \eqref{e306}
we have the claim.\\
\\
{\bf Theorem 8.2}(Main Theorem). {\it Let $|\omega |<3\,\beta ,$
$T>0$ and $u(x,\,t)$ be a solution of \eqref{e204} in the region
$\mathbb{R} \times [0,\,T]$ such that}
\begin{eqnarray}
\label{e807}u\in L^{\infty}([0,\,T]:\,H^{3}(W_{0\;L\;0}))
\end{eqnarray}
{\it for some $L\geq 2.$ Then }
\begin{eqnarray}
\label{e808}u\in L^{\infty}([0,\,T]:\,H^{3 + l}(W_{\sigma,\,L -
l,\,l})) \cap L^{2}([0,\,T]:\,H^{4 + l}(W_{\sigma,\,L - l - 1,\,l}))
\end{eqnarray}
{\it for all $0\leq l\leq L - 1$ and all $\sigma >0.$}\\
\\
{\it Remark.} If the assumption \eqref{e807} holds for all $L\geq
2,$ the solution is infinitely differentiable in the $x$-variable.
>From \eqref{e204} we have that the solution is $C^{\infty}$ in both
variables. We are also quantifying the gain of each derivative by
the degree of
vanishing of the initial data at infinity.\\
\\
{\it Proof.} We use induction on $\alpha.$ For $\alpha =3,$ let
$u$ be a solution of \eqref{e204} satisfying \eqref{e807}.
Therefore, $u_{t}\in L^{\infty}([0,\,T]:\,L^{2}(W_{0\;L\;0}))$
where $u\in L^{\infty}([0,\,T]:\,H^{3}(W_{0\;L\;0}))$ and
$u_{t}\in L^{\infty}([0,\,T]:\,L^{2}(W_{0\;L\;0})).$ Then $u\in
C([0,\,T]:\;L^{2}(W_{0\;L\;0}))\cap
C_{w}([0,\,T]:\,H^{3}(W_{0\;L\;0})).$ Hence, $u:[0,\,T]\longmapsto
H^{3}(W_{0\;L\;0})$ is a weakly continuous function. In
particular, $u(\,\cdot \,,\,t)\in H^{3}(W_{0\;L\;0})$ for all $t.$
Let $t_{0}\in (0,\,T)$ and $u(\,\cdot\,,\,t_{0})\in
H^{3}(W_{0\;L\;0}),$ then there are $\{u_{0}^{(n)}\}\subseteq
C_{0}^{\infty}(\mathbb{R})$ such that $u_{0}^{(n)}(\,\cdot
\,)\rightarrow u(\,\cdot\,,\,t_{0})$ in $H^{3}(W_{0\;L\;0}).$ Let
$u^{(n)}(x,\,t)$ be a unique solution of \eqref{e204} with
$u^{(n)}(x,\,t_{0})=u_{0}^{(n)}.$ Then by Theorem 5.1 and 5.2,
there exists $u$ in a time interval $[t_{0},\,t_{0} + \delta]$
where $\delta >0$ does not depend on $n$ and $u$ is a unique
solution of \eqref{e204}, $u^{(n)}\in L^{\infty}([t_{0},\,t_{0} +
\delta]:\,H^{3}(W_{0\;L\;0}))$ with $u^{(n)}(x,\,t_{0})\equiv
u_{0}^{(n)}(x)\rightarrow u(x,\,t_{0})\equiv u_{0}(x)$ in
$H^{3}(W_{0\;L\;0}).$ Now, by Theorem 7.1, we have
\begin{eqnarray*}
u^{(n)}\in L^{\infty}([t_{0},\,t_{0} +
\delta]:\,H^{3}(W_{0\;L\;0}))\cap L^{2}([t_{0},\,t_{0} +
\delta]:\,H^{4}(W_{\sigma ,\,L - 1,\,0}))
\end{eqnarray*}
with a bound that depends only on the norm of $u_{0}^{(n)}$ in
$H^{3}(W_{0\;L\;0}).$ Furthermore, Theorem 7.1 guarantees the
non-uniform bounds
\begin{eqnarray*}
\sup_{[t_{0},\,t_{0} + \delta]}\sup_{x}\,(1 + |x_{+}|)^{k}\,|\,
\partial^{\alpha}u^{(n)}(x,\,t)\,|<+\infty
\end{eqnarray*}
for each $n,\,k$ and $\alpha .$ The main inequality \eqref{e301}
and the estimate \eqref{e802} are therefore valid for each
$u^{(n)}$ in the interval $[t_{0},\,t_{0} + \delta].$ $\eta $ may
be chosen arbitrarily in its weight class \eqref{e801} and then
$\xi$ is defined by \eqref{e307} and the constant $c_{1},$
$c_{2},$ $c_{3},$ $c_{4}$ are independent of $n.$ From
\eqref{e301} and \eqref{e801} we have
\begin{eqnarray}
\label{e809}\sup_{[t_{0},\,t_{0} + \delta]}\int_{\mathbb{R}}\xi
\,|u_{\alpha}^{(n)}|^{2}dx + \int_{t_{0}}^{t_{0} +
\delta}\int_{\mathbb{R}}\eta \,|u_{\alpha + 1}^{(n)}|^{2}dx \leq c
\end{eqnarray}
where by \eqref{e802}, $c$ is independent of $n.$ The estimate
\eqref{e809} is proved by induction for $\alpha=3,\,4,\,5,\ldots $
Thus $u^{(n)}$ is also bounded in
\begin{eqnarray}
\label{e810}L^{\infty}([t_{0},\,t_{0} + \delta]:\,
H^{\alpha}(W_{\sigma,\,L - \alpha + 3,\,\alpha - 3})) \cap
L^{2}([t_{0},\,t_{0} + \delta]:\,H^{\alpha + 1}(W_{\sigma,\,L -
\alpha + 2,\,\alpha - 3}))
\end{eqnarray}
for $\alpha \geq 3.$ Since $u^{(n)}\rightarrow u$ in
$L^{\infty}([t_{0},\,t_{0} + \delta]:\,H^{3}(W_{0\;L\;0})).$ By
Corollary 5.3 it follows that $u$ belongs to the space \eqref{e810}.
Since $\delta $ is fixed, this result is valid over the whole
interval $[0,\,T].$

\end{document}